\documentclass[final,leqno,onefignum,onetabnum]{siamltex1213}

\usepackage{amsfonts}
\usepackage{amsmath,booktabs,ctable,threeparttable}
\usepackage{amssymb,amsfonts,boxedminipage}
\usepackage{algorithm}
\usepackage{algorithmic}

\newtheorem{remark}{Remark}[section]

\title{A Joint Bidiagonalization Based Algorithm for Large Scale
Linear Discrete Ill-posed
Problems in General-Form Regularization\thanks{This
work was supported in part by
the National Science Foundation of China (No. 11771249).}}

\author{Zhongxiao Jia and Yanfei Yang}
\author{Zhongxiao Jia\thanks{Department of Mathematical Sciences, Tsinghua University, 100084 Beijing,
China (\email{jiazx@tsinghua.edu.cn}).}
\and  Yanfei Yang\thanks{Department of Mathematical Sciences, Tsinghua University, 100084 Beijing,
China (\email{yangyf14@mails.tsinghua.edu.cn}).}}

\begin{document}
\maketitle
\slugger{sirev}{xxxx}{xx}{x}{x--x}

\begin{abstract}
Based on the joint bidiagonalization process of a large matrix pair $\{A,L\}$,
we propose and develop an iterative regularization algorithm for the large scale
linear discrete ill-posed problems in general-form regularization:
$\min\|Lx\| \ \mbox{{\rm subject to}} \ x\in\mathcal{S} =
\{x|\ \|Ax-b\|\leq \tau\|e\|\}$ with a Gaussian white noise $e$
and $\tau>1$ slightly, where $L$ is a regularization matrix.
Our algorithm is different from the hybrid one proposed by Kilmer {\em et al.},
which is based on the same process but solves
the general-form Tikhonov regularization problem:
$\min_x\left\{\|Ax-b\|^2+\lambda^2\|Lx\|^2\right\}$.
We prove that the iterates take the
form of attractive filtered generalized singular value decomposition (GSVD)
expansions, where the filters are given explicitly. This result and
the analysis on it show that the method must have the desired semi-convergence property
and get insight into the regularizing effects of the method.
We use the L-curve criterion or
the discrepancy principle to determine $k^*$.
The algorithm is simple and effective, and numerical experiments illustrate
that it often computes
more accurate regularized solutions than the hybrid one.
\end{abstract}



\pagestyle{myheadings}
\thispagestyle{plain}

\section{Introduction}

Consider the solution of the large scale linear discrete ill-posed problem
\begin{equation}\label{eq1}
\min_{x\in \mathbb{R}^{n}}\|Ax - b\| \quad {\rm or}
\quad Ax = b, \quad A \in \mathbb{R}^{m \times n}, \quad b \in \mathbb{R}^{m},
\end{equation}
where the norm $\|\cdot\|$ is the 2-norm of a vector,
the matrix $A$ is ill conditioned with its singular values
decaying to zero with no obvious gap between consecutive ones,
and the right-hand side $b = b_{true}+e$ is noisy and
assumed to be contaminated by a Gaussian white noise $e$, where $b_{true}$
is the noise-free right-hand side and $\|e\| < \|b_{true}\|$.
Such kind of problem arises in a variety of applications such as
computerized tomography, image deblurring,
signal processing, geophysics, heat
propagation, biomedical and optical imaging, groundwater modeling,
and many others; see, e.g.,
\cite{aster,berisha,engl93,hansen2007,kaipio,miller,nat1986}.
Since $b$ contains the noise $e$ and $A$ is extremely ill conditioned,
the naive solution $x_{naive} = A^{\dag}b$ is very large or huge in norm and
is a meaningless approximation to the true solution
$x_{true} = A^{\dag}b_{true}$,
where $\dag$ denotes the Moore-Penrose inverse of a matrix.
Therefore, one has to use regularization to obtain a best
possible approximation to $x_{true}$
\cite{hansen98,hansen10}.

Assume that $Ax_{true}=b_{true}$ and $m\geq n$. Then two essentially
equivalent dominating regularization approaches
are the following general-form regularization
\begin{align}\label{tik2}
\min\|Lx\| \ \ \ \mbox{{\rm subject to}} \ \ \ x\in\mathcal{S} =
\{x|\ \|Ax-b\|\leq \tau\|e\|\}
\end{align}
with some $\tau>1$  and the general-form Tikhonov regularization
\begin{align}\label{tik1}
\min_x\left\{\|Ax-b\|^2+\lambda^2\|Lx\|^2\right\},
\end{align}
where $L\in \mathbb{R}^{p\times n}$ is a regularization
matrix and $\lambda >0$ is the regularization parameter.
If $L=I_n$, the $n\times n$ identity matrix, then
\eqref{tik2} and
\eqref{tik1} are called standard-form
regularization problems.
If $A$ and $L$ satisfy
\begin{align}\label{nullcond}
\mathcal{N}(A)\cap\mathcal{N}(L)=\{\mathbf{0}\},\ {\rm i.e.},\
{\rm rank}\left(
\begin{array}{c}
A \\
L \\
\end{array}
\right)=n,
\end{align}
the solution to \eqref{tik1} is unique, where $\mathbf{0}$ denotes
the zero vector of dimension $n$. In practical applications,
$L$ is typically chosen as a scaled approximation of the first or
second derivative operator \cite{hansen98,hansen10}.

For $L\not=I_n$, \eqref{tik2} and \eqref{tik1} can be transformed to their
standard forms with $L=I_n$ and $A$ replaced by $AL_A^{\dagger}$,
where
$$
L_A^{\dagger}=(I-(A(I-L^{\dagger}L))^{\dagger}A)L^{\dagger}
$$
is {\em $A$-weighted pseudoinverse of $L$} and
$L_A^{\dagger}=L^{\dagger}$ when $p\geq n$;
see \cite{hansen98} for details.
This is computationally viable and attractive
if not much effort is needed by applying $L_A^{\dagger}$, e.g.,
when $L$ is banded with small bandwidth and has a known null space; 
we refer the reader to, e.g., \cite{chung2017,chung2018,gazzola}
for some available algorithms and codes.
In many practical applications, however,
such transformation is computationally unfeasible.
This is often the case for two or three dimensional case, e.g.,
where $L$ has no special structure
or is the sum of Kronecker products such that $L_A^{\dagger}$
is expensive to use.

There have been some available randomized algorithms and Krylov
subspace type methods for
solving the regularization problem \eqref{tik2} or \eqref{tik1}
when the application of $L_A^{\dagger}$ is computationally unfeasible.
For example, Jia and Yang \cite{jiayang18} have proposed and developed
efficient randomized SVD algorithms for solving \eqref{tik2} effectively, where
a large scale least squares problem is iteratively solved with low or modest
accuracy at each step.
Several Krylov subspace type methods have been presented to solve \eqref{tik1}.
To shed light on a few existing Krylov subspace type methods,
it is important and necessary to keep in mind the basic and core requirements
for a regularization method:
As far as solving \eqref{tik1} is concerned, (i) a good regularized solution
must capture the dominant generalized singular value decomposition
(GSVD) components of the matrix pair $\{A,L\}$
and meanwhile suppress those corresponding to small generalized singular values,
and (ii) the generalized right singular vectors of the matrix pair $\{A,L\}$
form a more suitable basis to express a regularized solution.
See \cite{hansen98,hansen10,kilmer2007} for details.

Based on some generalized Arnoldi process proposed first
by Li and Ye  \cite{li2003} for the solution
of quadratic eigenvalue problems, Reichel \emph{et al.} \cite{reichel12}
present a hybrid iterative algorithm for solving \eqref{tik1}.
The generalized Arnoldi process successively reduces
the square matrix pair $\{A,L\}$ to a sequence of small
matrix pairs $\{H_A,H_L\}$, where the
projection matrix $H_L$ quickly becomes
full as the number of iterations $k$ increases. The reduction exploits
only $A$ and $L$ but does not use their transposes, so the information on
$A^T$ and $L^T$ is lacking. This may make the underlying subspace unable
to capture a dominant generalized right singular subspace. In fact,
for $L=I_n$, Hansen in his book \cite[p.126]{hansen10}
insightfully summarizes that Arnoldi process based methods, such as
RRGMRES, mix the SVD components in each iteration, their success is
highly problem dependent, and they can successful when the mixing of the SVD
components is weak, e.g.,
$A$ is (nearly) symmetric. It follows from the above
that generalized Arnoldi methods have similar limitations. Indeed,
we have observed from \cite{reichel12} that
the regularized solutions behaved quite irregularly
when $A$ is nonsymmetric; see Example 5.1 there. This makes it
very hard to stop the algorithm properly.

Hochstenbach \emph{et al.} \cite{hoch15} propose an extended Golub-Kahan
bidiagonalization process to reduce the matrix pair $\{A,L\}$
to a sequence of small matrix pairs $\{H_{k+1,k},K_{k,k}\}$ with $H_{k+1,k}$
and $K_{k,k}$ being upper Hessenberg and triangular, respectively.
The process degenerates to the standard Golub-Kahan bidiagonalization
process when $L=I_n$.
They then develop a hybrid projection algorithm for solving \eqref{tik1},
which is projected onto a sequence of generalized Krylov subspaces
generated by $A^Tb$ and the matrices $A^TA$ and $L^TL$
simultaneously. The underlying solution subspace contains
the much lower dimensional standard Krylov subspace generated by $A^Tb$ and
$A^TA$ and that generated by $A^Tb$ and $L^TL$
and meanwhile includes too many possibly useless vectors.
Precisely, at iteration $k$, the process produces
only $\lfloor \log_2 k\rfloor$ dimensional Krylov subspaces generated
by $A^Tb$ and $A^TA$ and by $A^Tb$ and $L^TL$, respectively,
and each iteration
computes the matrix-vector products with $A^T, A, L^T$ and $L$ and uses
{\em longer} recurrences during orthonormalization of basis vectors
as the iteration {\em proceeds};
see \cite{hoch15,zwaan}. In other words, the dimension of
generalized Krylov subspace increases {\em exponentially}
when generating a standard Krylov subspace generated by each of
$A^TA$ and $L^TL$, respectively.
For example, if only fifteen dimensional Krylov subspaces
generated by $A^TA$ and $L^TL$ are needed, one has to perform
the extended Golub-Kahan bidiagonalization
process $2^{15}= 32768$ steps if
$2^{15}\leq \min\{n,p\}$, which is very expensive. In
the meantime, we should notice that the Krylov subspace
generated by $A^Tb$ and $A^TA$ favors dominant right
singular vectors of the single $A$, which may bear no relation
to those {\em desired} dominant generalized ones of the matrix pair
$\{A,L\}$. In the meantime, it is not yet clear what the
whole generalized Krylov subspace favors.
We mention that, except Example 4.5 in \cite{hoch15}
where it is unclear whether or not $A$ is symmetric,
all the other test matrices $A$ in \cite{hoch15,zwaan} are symmetric,
in which case the method is more possible to succeed.


Zha \cite{zha1996} proposes a joint bidiagonalization process that
successively reduces the matrix pair $\{A,L\}$ to upper
bidiagonal forms.
Based on the process, Zha proposes a joint bidiagonalization method
for computing a few largest or smallest generalized singular values and the
associated singular vectors of a large matrix pair $\{A,L\}$.

Kilmer \emph{et al.} \cite{kilmer2007} adapt Zha's joint bidiagonalization
to discrete linear ill-posed problems in
general-form regularization and develop a joint bidiagonalization
process that successively reduces the matrix pair $\{A,L\}$ to lower and upper
bidiagonal forms. Based on the process, they propose a hybrid projection
method for solving \eqref{tik1}.
It is argued in \cite{kilmer2007} that the underlying solution subspaces are
{\em legitimate} since they appear to be more directly related to
the generalized right singular vectors of $\{A,L\}$. Unlike the extended
Golub-Kahan bidiagonalization process, the $k$-step joint bidiagobalization
process does not make any possible waste, but at each step
it needs the solution of a large scale linear least squares problem with
the coefficient matrix $(A^T,L^T)^T$ that is supposed to be solved
iteratively, called inner iteration. Therefore, joint bidiagonalization
forms an inner-outer iterative process. Fortunately,
$(A^T,L^T)^T$ is typically well conditioned as $L$ is typically so
in applications \cite{hansen98,hansen10}. In these cases,
the LSQR algorithm \cite{lsqr1982} can
efficiently solve the mentioned least squares problems. At each iteration of
the hybrid projection method \cite{kilmer2007},
one solves a small projected general-form Tikhonov regularization problem.
Finally, one solves a large scale least squares problem
with the coefficient matrix $(A^T,L^T)^T$
to form a regularized solution. The outer iteration proceeds until
the regularized solutions stagnate, that is, their accuracy
cannot be improved. In the hybrid method, the iteration
number do not play the role of the regularization parameter, and
one needs to determine an optimal Tikhonov regularization parameter
for a projected general-form regularization problem generated
at each outer iteration.

It is well known from, e.g., \cite{hansen98,hansen10},
that any regularization is based on an underlying requirement
that the discrete Picard condition for a given problem is satisfied,
only under which can one compute a useful regularized solution with some accuracy.
Notice that a small projected general-form  Tikhonov regularization
problem is solved at each outer iteration \cite{kilmer2007},
where one needs to determine an optimal regularization parameter for
the small problem itself. This is also the case for any hybrid projection
method for solving \eqref{tik1}. One fundamental fact that is crucial but
has received little attention until the work
\cite{gazzola16,renuat17} is: With $L=I_n$, that \eqref{tik1} satisfies
the discrete Picard condition does not mean that the
projected problems fulfill discrete Picard conditions too,
and known sufficient conditions for the projected problems to satisfy
the discrete Picard conditions require that the singular values of
the matrices involved in the projected problems
approximate the large singular values of $A$ in natural order.
An adaption of this result to $L\not=I_n$
says that the projected problems are guaranteed
to inherit the discrete Picard conditions when the generalized singular values
of the projected matrix pairs approximate the large generalized
singular values of the matrix pair $\{A,L\}$ in natural order.
However, it has been proved in \cite{jia16,jiarank}
that, for $L=I_n$ and LSQR, the approximations in natural
order can be guaranteed only for severely ill-posed problems and
some moderately ill-posed ones and such property generally
does not hold when the singular values of $A$ or the generalized ones
of $\{A,L\}$ decay slowly, e.g., those matrices in
mildly ill-posed problems. For the definition
of severely, moderately and mildly ill-posed problems, see
\cite{hansen98,hansen10}.

If the above sufficient conditions are met, the iterative algorithm
used resembles the truncated GSVD (TGSVD) method \cite{hansen98}
until the occurrence of semi-convergence.
At this time, a best regularized solution has been already found and is
as accurate as the best TGSVD solution, which is a best
regularized solution of \eqref{eq1} in the sense of
the regularization formulation \eqref{tik2}. Therefore, one only needs
to determine the semi-convergence point by some parameter-choice
methods, e.g., the L-curve criterion and the discrepancy principle.
We should be aware that, for a hybrid iterative method,
whenever the discrete Picard conditions for the projected
problems fail to satisfy
or are satisfied poorly, i.e., their solution norms (very) large,
optimal regularization parameters for them
are poorly defined. In this case, a direct consequence is that
the regularized solutions may exhibit irregular behavior.
When developing a hybrid LSQR variant,
by requiring that the singular values of the projected matrices
approximate the large singular values of $A$ in natural order,
Renaut {\em et al.} \cite{renuat17} prove that an optimal regularization
parameter for each projected problem can be reliably determined by a weighted
generalized cross-validation (WGCV) parameter-choice method and it converges
to the global optimal regularization parameter $\lambda_{opt}$ for \eqref{tik1}
as outer iterations proceed, so that the regularized solutions
ultimately stagnate.

In this paper, based on joint bidiagonalization process \cite{kilmer2007},
instead of developing any hybrid projection based algorithm,
we will propose a pure projection iterative algorithm
for solving \eqref{tik2} other than the equivalent
\eqref{tik1} as done in \cite{kilmer2007}. Our algorithm is much simpler
than the hybrid one \cite{kilmer2007}, and the iteration number $k$
plays the role of the regularization parameter.
First, we make use of the joint bidiagonalization process to
project \eqref{tik2} onto a sequence of low dimensional subspaces and
obtain a sequence of projected problems, which involve
matrix pairs of small size. Then at each outer iteration we solve a projected
problem. Remarkably, we find out that the solution of each of them
reduces to an {\em ordinary} small least squares problem with the coefficient
matrix being single lower bidiagonal other than the matrix pair,
so that it is very cheap and reliable to solve it
by the QR factorization at $\mathcal{O}(k)$ flops.
Another big benefit is that we no longer determine
the optimal regularization parameter for each projected problem, which itself
may be poorly defined. We make a theoretical analysis on the proposed method
and establish a number of results. Most importantly,
we first prove that the iterates obtained by our algorithm take the
form of filtered GSVD expansions and are expressed {\em explicitly}
in the generalized right singular vector basis of $\{A,L\}$,
a desired and insightful property. Then we analyze the filters in the expansions,
shed light on the regularizing effects of the algorithm, and prove that
the iterates capture dominant GSVD components of $\{A,L\}$, as desired.
The results indicate that our
algorithm must have the typical semi-convergence property:
As the joint bidiagonalization process proceeds, more and more
dominant generalized singular components of $\{A,L\}$ are captured, and
the regularized solutions converge to the true solution
$x_{true}$ of \eqref{eq1} until some iteration, after which
the regularized solutions start to be deteriorated by the noise $e$ and
instead converge to $x_{naive}$.

As it will turn out, since the residual norms monotonically decrease
and the semi-norms of solutions practically increase monotonically,
we can use the L-curve criterion and the discrepancy
principle to estimate
the optimal regularization parameter $k^*$, at which the
semi-convergence occurs.

We will numerically compare our algorithm with the hybrid
algorithm in \cite{kilmer2007}, in which we make use of the GCV
and WGCV parameter choice methods to determine
an optimal regularization parameter
for each small projected problem. We are primarily
concerned with the accuracy of the best regularized solutions by our algorithm
and the hybrid one. The experiments on several real-world problems
will illustrate the superiority of our algorithm.

The paper is organized as follows. In Section 2, we
overview GSVD, the general-form Tikhonov
regularization method and the TGSVD method, and present
the joint bidigonalization process of $\{A,L\}$. In Section 3, we
describe the hybrid method in \cite{kilmer2007}.
In Section 4, we propose our joint
bidiagonalization based method for solving \eqref{tik2}. In Section 5, we
make a theoretical analysis on it.
In Section 6, we consider the practical determination of
the optimal regularization parameter. Numerical experiments are presented in
Section 7. Finally, we conclude the paper in Section 8.

\section{GSVD, regularization methods and joint bidiagonalization}

In this section, we provide some necessary background. We describe GSVD, the TGSVD
method, the filtered GSVD method,
and the joint bidiagonalization process proposed in \cite{zha1996}
and developed in \cite{kilmer2007}.

Consider the compact QR factorization
\begin{align}\label{qrfact}
\left(
  \begin{array}{c}
    A \\
    L \\
  \end{array}
\right)=QR,
\end{align}
where $Q=\left(
           \begin{array}{c}
             Q_A \\
             Q_L \\
           \end{array}
         \right)
\in\mathbb{R}^{(m+p)\times n}$ is column orthonormal with
$\ Q_A\in\mathbb{R}^{m\times n},\ Q_L\in\mathbb{R}^{p\times n}$,
and $R\in\mathbb{R}^{n\times n}$ is upper triangular and nonsingular
because of the assumption \eqref{nullcond}.
We have $A=Q_AR, \ L = Q_LR$, and $Q_A^TQ_A + Q_L^TQ_L = I_n$.

Let the CS decomposition of the matrix pair $\{Q_A, Q_L\}$ be
\begin{align}\label{csd}
Q_A = P_ACW^T,\ \ Q_L=P_LSW^T,
\end{align}
where $P_A\in\mathbb{R}^{m\times m}$, $P_L \in\mathbb{R}^{p\times p}$,
and $W\in\mathbb{R}^{n\times n}$ are orthogonal, and
$C\in\mathbb{R}^{m\times n}$ and $S\in\mathbb{R}^{p\times n}$ are
diagonal matrices satisfying $C^TC+S^TS=I_n$; see \cite[Section 4.2]{bjorck96}.
Then the GSVD of $\{A,L\}$ is
\begin{align}\label{gsvd}
A = P_ACG^{-1},\ \ L=P_LSG^{-1}
\end{align}
with $G=(g_1,g_2,\ldots,g_n)=R^{-1}W\in \mathbb{R}^{n\times n}$, and the vectors
$g_i$ are the generalized right singular vectors of $\{A,L\}$.
Following the unconventional but more convenient way \cite{kilmer2007},
we order the entries of the diagonal matrices $C$ and $S$ so that
\begin{equation}\label{svda}
1\geq c_1\geq\cdots \geq c_{\min\{n,p\}}\geq 0, c_{\min\{n,p\}+1}=\cdots =c_n=1,
\end{equation}
\begin{equation}\label{svdl}
0 \leq s_1 \leq\cdots \leq s_{\min\{n,p\}} \leq 1.
\end{equation}
By the GSVD \eqref{gsvd}, the general-form Tikhonov solution $x_{\lambda}$
to \eqref{tik1} takes a filtered GSVD expansion:
\begin{align}\label{tiksol}
x_{\lambda}& = (A^TA + \lambda^2 L^TL)^{-1}A^T b = G(C^T
C + \lambda^2 S^TS)^{-1}C^TP_A^Tb \notag\\
&=\sum_{i=1}^{\min\{n,p\}}\frac{c_i^2}{c_i^2+\lambda^2 s_i^2}
\frac{p_{i,A}^Tb}{c_i}g_i+\sum_{i=\min\{n,p\}+1}^n p_{i,A}^T b g_i,
\end{align}
where $P_A=(p_{1,A},p_{2,A},\ldots,p_{m,A})$, $f_i
=\frac{c_i^2}{c_i^2+\lambda^2 s_i^2}$
are filters, and the second term lies in the null space $\mathcal{N}(L)$
of $L$, which is spanned by the vectors $g_{\min\{n,p\}+1},\ldots,g_n$.
We address that the regularization does not affect the second term.
This is simply the filtered GSVD method for solving \eqref{tik1}.

The discrete Picard condition \cite{hansen98} states that
the Fourier coefficients $|p^T_{i,A} b|$ must, on average, decay faster than
the $c_i$. Hence the $|p_{i,A}^Tb|$ decay until
the Gaussian white noise $e$ dominates the $|p_{i,A}^Tb|$,
that is, the $|p_{i,A}^Tb|\approx
|p_{i,A}^T e|$ stagnates and is dominated by $e$
after $i>k_0$ for some $k_0$, while
$|p_{i,A}^Tb|\approx |p_{i,A}^Tb_{true}|$ is dominated by
$b_{true}$, $i=1,2,\ldots,k_0$, where $k_0$
is called the transition or cutting-off point.
Therefore, a good regularized solution $x_{\lambda}$ must capture
the $k_0$ dominant GSVD components of $\{A,L\}$ and
meanwhile dampen those for $i>k_0$ as much as possible.
An optimal regularization parameter $\lambda_{opt}$
can be determined by some parameter-choice methods, e.g.,
the discrepancy principle, the L-curve criterion, and the generalized
cross validation (GCV) or weighted GCV (WGCV) method;
see \cite{hansen98,hansen10} and also \cite{berisha,hansen2007}.

Alternatively, making use of the GSVD of $\{A,L\}$, one can develop
the TGSVD method and computes a sequence of the TGSVD solutions
\begin{equation}\label{tgsvdsol}
x_k^{tgsvd}=
\sum_{i=1}^k\frac{p_{i,A}^Tb}{c_i}g_i+
\sum_{i=\min\{n,p\}+1}^n p_{i,A}^T b g_i,\ \ k=1,2,\ldots,\min\{n,p\},
\end{equation}
where the first term consists of the first $k$ dominant GSVD
components of $\{A,L\}$. The TGSVD solution takes a special filtered GSVD
expansion, where the filters $f_i=1$ for $i=1,2,\ldots,k$ and $f_i=0$ for
$i=k+1,\ldots,\min\{n,p\}$.
Under the discrete Picard condition, the TGSVD method exhibits
semi-convergence: $x_k^{tgsvd}$ and $Lx_k^{tgsvd}$ converge to $x_{true}$ and
$Lx_{true}$ for $k \leq k_0$, afterwards they diverge
and instead converge to $x_{naive}$ and $Lx_{naive}$, respectively.
A best possible TGSVD solution $x_{k_0}^{tgsvd}$ is thus obtained for $k=k_0$.

We notice that the second terms in \eqref{tiksol} and \eqref{tgsvdsol} are
the same and they disappear when $p\geq n$. For later use, we write them as
\begin{equation}\label{second}
g_{\perp}=\sum_{i=\min\{n,p\}+1}^n p_{i,A}^T b g_i \in \mathcal{N}(L).
\end{equation}

Now we review a procedure that jointly
diagonalizes the matrix pair $\{A,L\}$ to lower and upper
bidiagonal forms. Applying the BIDIAG-1 algorithm
and BIDIAG-2 algorithm in \cite{lsqr1982} to $Q_A$ and $Q_L$, respectively,
which are the lower and upper Lanczos bidiagonalization processes,
we can reduce $Q_A$ and
$Q_L$ to lower and upper bidiagonal forms, respectively. The two processes
can be written in matrix form:
\begin{align}
&Q_AV_k=U_{k+1}B_k,\ \ Q_A^TU_{k+1}=V_kB_k^T+\alpha_{k+1}v_{k+1}e_{k+1}^T,\label{diag1}\\
&Q_L\widehat{V}_k=\widehat{U}_k\widehat{B}_k,\ \ Q_L^T\widehat{U}_k
=\widehat{V}_k\widehat{B}_k^T+\widehat{\beta}_k\widehat{v}_{k+1}e_k^T,\label{process}
\end{align}
where $e_{k+1}$ and $e_k$ are the $(k+1)$th and $k$th canonical vectors of
dimensions $k+1$ and $k$, respectively,
\begin{equation}\label{uk22}
U_{k+1}=(u_1,\ldots, u_{k+1})\in\mathbb{R}^{m\times (k+1)},
\ \ \widehat{U}_k=(\widehat{u}_1,\ldots\widehat{u}_k)\in\mathbb{R}^{p\times k},
\end{equation}
and
\begin{equation}\label{vk22}
V_k=(v_1,\ldots,v_k)\in\mathbb{R}^{n\times k}, \ \
\widehat{V}_k=(\widehat{v}_1,\ldots,\widehat{v}_k)\in\mathbb{R}^{n\times k}
\end{equation}
are column orthonormal, and
\begin{equation}\label{bb}
B_k=\left(
        \begin{array}{cccc}
          \alpha_1 &  &  &  \\
          \beta_2  & \alpha_2 &  &  \\
                  &  \beta_3 & \ddots &  \\
                  &      & \ddots & \alpha_k \\
           &  &  & \beta_{k+1} \\
        \end{array}
      \right)\in\mathbb{R}^{(k+1)\times k }, \ \
\widehat{B}_k=\left(
        \begin{array}{cccc}
          \widehat{\alpha}_1 & \widehat{\beta}_1  &  &  \\
           & \widehat{\alpha}_2 &  \ddots&  \\
                  &   & \ddots &   \widehat{\beta}_{k-1}\\
                  &      & & \widehat{\alpha}_k \\
                  \end{array}\right)\in\mathbb{R}^{k\times k }
\end{equation}
are lower bidiangonal and upper bidiagonal, respectively.
Zha \cite{zha1996} and Kilmer \emph{et al.} \cite{kilmer2007} have
investigated the relationships
between $V_k$ and $\widehat{V}_k$ defined in \eqref{vk22} and between
$B_k$ and $\widehat{B}_k$ defined in \eqref{bb}, respectively, and
they have established the following result.

\begin{theorem}\label{zha}
If $v_1=\widehat{v}_1$ in \eqref{vk22}, then
\begin{equation}
\widehat{v}_{j+1}=(-1)^jv_{j+1}, \ \
\widehat{\alpha}_j\widehat{\beta}_j=\alpha_{j+1}\beta_{j+1},\ \
j=1,\ldots,k.
\end{equation}
\end{theorem}

A combination of \eqref{diag1}--\eqref{bb}, Theorem \ref{zha} and
the QR factorization \eqref{qrfact} shows
that $A$ and $L$ can be jointly bidiagonalized \cite{kilmer2007,zha1996},
as summarized below.

\begin{theorem}
Assume that $A \in \mathbb{R}^{m \times n}$ and $L \in \mathbb{R}^{p \times n}$
with $m\geq n$. Then there exist orthogonal matrices $U \in \mathbb{R}^{m \times m}$,
$\widehat{U}\in \mathbb{R}^{p \times p}$ and $V \in \mathbb{R}^{n \times n}$,
and a lower bidiagonal $B \in \mathbb{R}^{m \times n}$, an upper bidiagonal $\bar{B}
\in \mathbb{R}^{p \times n}$, and an invertible $Z$ such that
\begin{align}\label{bid0}
A&= UBZ^{-1},\\
L&= \widehat{U}\bar{B}Z^{-1},\label{bid01}
\end{align}
where $Z = R^{-1}V$ and $\bar{B} = \widehat{B}D$ with $D=\diag(1,-1,1,-1,\ldots)$,
and the remaining
matrices are obtained by running joint bidiagonalization to completion.
In particular, when $p< n$, the columns $p+1,\ldots, n$ of $\bar{B}$
contain only zero entries.
\end{theorem}

From \eqref{bid0} and \eqref{bid01}, we obtain $k$-step joint bidiagonalization
relations
\begin{align}\label{bid2}
AZ_k &= U_{k+1}B_k,\\
LZ_k &= \widehat{U}_k\bar{B}_k\label{bid22},
\end{align}
where $Z_k\in\mathbb{R}^{n\times k}$ is the first $k$ columns of $Z$,
and $B_k$ and $\bar{B}_k$ are the first $(k+1)\times k$ and $k\times k$
submatrices of $B$ and $\bar{B}$, respectively.

For $A$ and $L$ large, the QR factorization \eqref{qrfact} is impractical.
In order to avoid explicitly computing $Q_A$ and $Q_L$, inspired by Zha's
work \cite{zha1996}, Kilmer {\em et al.} \cite{kilmer2007}
develop a joint bidiagonalization (JBD) process,
denoted by Algorithm~\ref{alg:1},
to compute the matrices in \eqref{uk22}--\eqref{bb},
in which $0_p$ denotes the zero vector of dimension $p$.

\begin{algorithm}[!htp]
\caption{\ $k$-step joint bidiagonalization (JBD) process.}
\begin{algorithmic}[1]\label{alg:1}
\STATE $\beta_1u_1 = b$, $\beta_1=\|b\|$.\\
\STATE $\alpha_1\widetilde{v}_1=QQ^T
\left(\begin{array}{c}
                     u_1 \\
                     0_p \\
                   \end{array}
                 \right).$\\
\STATE $\widehat{\alpha}_1\widehat{u}_1=\widetilde{v}_1(m+1:m+p)$\\
\FOR{$i=1, 2,\ldots,k$}
\STATE $\beta_{i+1}u_{i+1}=\widetilde{v}_i(1:m)-\alpha_iu_i$.\\
\STATE $\alpha_{i+1}\widetilde{v}_{i+1}=QQ^T\left(
                   \begin{array}{c}
                     u_{i+1} \\
                     0_p \\
                   \end{array}
                 \right)-\beta_{i+1}\widetilde{v}_i$.\\
\STATE $\widehat{\beta}_i=(\alpha_{i+1}\beta_{i+1})/\widehat{\alpha}_i$.\\
\STATE $\widehat{\alpha}_{i+1}\widehat{u}_{i+1}
         =(-1)^i\widetilde{v}_{i+1}(m+1:m+p)-\widehat{\beta}_i\widehat{u}_i$.
\ENDFOR
\end{algorithmic}
\end{algorithm}

Let $\widetilde{u}_i=\left(\begin{array}{c} u_i \\ 0_p \\\end{array} \right)$.
At each iteration $i=1,2,\ldots,k+1$, Algorithm~\ref{alg:1} needs to compute
$QQ^T\widetilde{u}_i$, which is not accessible since $Q$ is not available.
However, notice that $QQ^T\widetilde{u}_i$ is nothing but the orthogonal projection
of $\widetilde{u}_i$ onto the column space of $\left(\begin{array}{c} A \\ L \\
\end{array} \right)$, which means that $QQ^T\widetilde{u}_i=\left(\begin{array}{c} A
\\ L \\\end{array} \right)\widetilde{x}_i$. Clearly, $\widetilde{x}_i$
is the solution to the least squares problem:
\begin{equation}\label{innerleast}
\widetilde{x}_i=\arg\min_{\widetilde{x}\in \mathbb{R}^n}
\left\|\left(\begin{array}{c} A \\ L \\\end{array} \right)
\widetilde{x}-\widetilde{u}_i\right\|.
\end{equation}
Since the least squares problem is large scale,
it is generally only feasible to solve it by an iterative solver, e.g.,
the most commonly used LSQR algorithm \cite{lsqr1982}.
Here we have two remarks.
\begin{remark}
Theoretically, at outer iteration $k$,
the inner least square problem \eqref{innerleast}
is solved accurately in order to guarantee that
\eqref{bid2} and \eqref{bid22} hold exactly. It is unknown
whether or not the solution accuracy
can be relaxed by allowing possibly large inexactness
in the algorithm \cite{kilmer2007} and ours to be presented later.
This issue is certainly complicated.  We do not investigate it
in the current paper. In finite precision arithmetic,
we suppose that \eqref{innerleast} is solved by the Matlab function
{\sf lsqr.m} with the default stopping criterion $10^{-6}$.
\end{remark}

\begin{remark}
To ensure the numerical orthogonality of the computed $U_{k+1},\
\widehat{U}_k$ and $V_k$, we use one step reorthogonalization during the process
in implementations.
\end{remark}

\section{The hybrid projection based method in \cite{kilmer2007}}

Algorithm~\ref{alg:1} takes
\begin{align}\label{bid1}
U_{k+1}(\beta_1e_1)&=b,\ \ \ \beta_1=\|b\|,
\end{align}
where $e_1$ is the first canonical vector of dimension $k+1$.
For a given regularization parameter $\lambda$, the hybrid projection
based method in \cite{kilmer2007} seeks the solution $x_k^{\lambda}\in
{\rm span}\{Z_k\}$ such that
$$
\min_{x\in {\rm span}\{Z_k\}}\{\|Ax-b\|^2+\lambda^2\|Lx\|^2\}=
\|Ax_k^{\lambda}-b\|^2+\lambda^2\|Lx_k^{\lambda}\|^2.
$$
Exploit \eqref{bid1}, \eqref{bid2} and \eqref{bid22}, and write
\begin{equation}\label{xk}
x_k^{\lambda}=Z_ky_k^{\lambda}.
\end{equation}
It is direct to justify that
\begin{align*}
\|Ax_k^{\lambda}-b\|^2+\lambda^2\|Lx_k^{\lambda}\|^2
&=\|B_ky_k^{\lambda}-\beta_1e_1\|^2+\lambda^2\|\bar{B}_ky_k^{\lambda}\|^2 \notag \\
&=\min_{y}\{\|B_ky-\beta_1e_1\|^2+\lambda^2\|\bar{B}_ky\|^2\}.
\end{align*}
Therefore, at iteration $k$ the hybrid projection
based method in \cite{kilmer2007} solves a projected general-form Tikhonov
regularization problem
\begin{align}\label{pro1}
\min_{y}\{\|B_k y-\beta_1e_1\|^2+\mu_k^2\|\bar{B}_k y\|^2\},
\end{align}
where the new notation $\mu_k>0$ is introduced to specialize
the regularization parameter for the projected problem at
iteration $k$. The key is the
determination of an optimal regularization $\mu_{kopt}$ for \eqref{pro1}.
Following the results of Renaut \emph{et al.} \cite{renuat17}
with $L=I_n$ to \eqref{pro1},
the optimal $\mu_{kopt}$ determined by the GCV or WGCV method converges to
the global optimal regularization parameter $\lambda_{opt}$ for \eqref{tik1}
as $k$ increases under the assumption that the generalized singular values
of $\{B_k,\bar{B}_k\}$
approximate the large singular values of $\{A,L\}$ in natural order; if the
assumption is not fulfilled, such convergence may fail, implying
that the regularized solution
$x_k^{\mu_{kopt}}$ may behave irregularly and does not stabilize for $k$
sufficiently large. As a consequence, it may be hard to stop the hybrid algorithm
properly, and even for $k$ sufficiently large the regularized solution
$x_k^{\mu_{kopt}}$ may not be as accurate
as $x_{\lambda_{opt}}$, the best regularized solution to \eqref{tik1}
associated with $\lambda=\lambda_{opt}$.

Now we show how to compute $x_k^{\lambda}$ when $y_k^{\lambda}$ is known.
Let $\widetilde{V}_k=(\widetilde{v}_1,\ldots, \widetilde{v}_k)
\in\mathbb{R}^{(m+p)\times k}$
be generated by Algorithm \ref{alg:1}. Then
\begin{equation}\label{vk}
\widetilde{V}_k=QV_k=QR(R^{-1}V_k)
=\left(\begin{array}{c} A \\L \\\end{array}\right)Z_k,
\end{equation}
from which, \eqref{qrfact} and \eqref{xk} it follows
that
\begin{equation}\label{inner}
\left(\begin{array}{c} A \\L \\\end{array}\right)
x_k^{\lambda}=\widetilde{V}_ky_k^{\lambda}.
\end{equation}
Kilmer {\em et al.} \cite{kilmer2007} show
that one only needs to form $x_k^{\lambda}$ explicitly when it
is accepted as the final regularized solution.

Regarding the determination of $\lambda_{opt}$,
other than determining $\mu_{kopt}$ for each projected problem \eqref{pro1},
Kilmer {\em et al.} \cite{kilmer2007} use the L-curve criterion
to tentatively estimate $\lambda_{opt}$.
They assume that a set of $\lambda$-values is prescribed,
and derive some update formulas for all the quantities,
including regularized solutions, residual norms, and the semi-norms
$\|Lx_k^{\lambda}\|$, which can be efficiently computed for the a-prior
set of $\lambda$-values. For sufficiently large $k$ at which
all the needed dominant GSVD components of $\{A,L\}$ are
thought to have been captured, drawing the picture of
$(\log \|Ax_k^{\lambda}-b\|,\log \|Lx_k^{\lambda}\|)$ for
the given set of $\lambda$-values,
they attempt to obtain a L-curve and pick up the $\lambda$-value at the corner
as an approximation $\lambda_{opt}$.

Their approach to determining $\lambda_{opt}$
faces two challenging issues: The first is how to
effectively determine a sufficiently large $k$, and the second is how
to choose a good a-prior set of $\lambda$-values which include
the optimal regularization parameter $\lambda_{opt}$ or its
good approximation. As a matter of fact, the first issue is common in
any hybrid LSQR algorithm, and there has been no reliable approach to determine
it. In our implementation, we will take a regular manner, as done in, e.g.,
\cite{hansen98,renuat17}, and determine $\mu_{kopt}$
for each \eqref{pro1} by the GCV code \cite{hansen2007} and the
WGCV code adapted from \cite{berisha}, both of which need to
compute the GSVD of $\{B_k,\bar{B}_k\}$ at cost of $\mathcal{O}(k^3)$ flops.

\section{Our joint bidiagonalization based algorithm}

Instead of solving \eqref{tik1}, we now present a joint
bidiagonalization based algorithm for solving \eqref{tik2}, which is a
purely iterative regularization method {\em without} explicitly
regularizing projected problems at each outer iteration. In the algorithm,
the iteration number $k$ plays the role of the regularization
parameter. We will establish a number of important results, which get
insight into the regularizing effects of the proposed method
and particularly prove that the method has the desired semi-convergence.

Still, we seek $x_k\in {\rm span}\{Z_k\}$ and write it in the form
\begin{equation} \label{xkyk}
x_k=Z_k y_k.
\end{equation}
We aim to project the original large scale regularization problem \eqref{tik2}
onto a sequence of low dimensional subspaces ${\rm span}\{Z_k\}$
and compute regularized solutions $x_k$ from them.
We will consider when to terminate and how to practically determine a best
regularized solution in Section 6.

From \eqref{xkyk}, the above projection is equivalent
to replacing $A$ and $L$ by $AZ_k$ and $LZ_k$
in \eqref{tik2} and solving the {\em reduced} general-form regularization problem
\begin{align}\label{reduced}
\min\|LZ_k y\| \ \ \ \mbox{{\rm subject to}} \ \ \ y\in
\{y|\ \|AZ_k y-b\|=\min\}
\end{align}
for $y_k$, starting with $k=1$ onwards. Make use of \eqref{bid1}, \eqref{bid2}
and \eqref{bid22}. Then \eqref{reduced} becomes
the {\em reduced} general-form regularization problem
\begin{align}\label{pro2}
\min\|\bar{B}_ky \| \ \ \ \mbox{{\rm subject to}}
\ \ \ y\in \left\{y |\ \|B_ky -\beta_1e_1\|=\min\right\}
\end{align}
starting with $k=1$ onwards. After the solution $y_k$ for \eqref{pro2} is computed,
in terms of \eqref{vk}, \eqref{qrfact} and $x_k=Z_k y_k$,
we then compute $x_k$ by solving
\begin{equation}\label{xk2}
\left(\begin{array}{c} A \\L \\\end{array}\right)x_k=\widetilde{V}_ky_k.
\end{equation}

Now let us investigate the solution of the constrained problem \eqref{pro2}.
As it will turn out, \eqref{pro2} amounts to
an ordinary unconstrained linear least squares problem, as shown below.

\begin{theorem}\label{equal}
Assume that Algorithm~\ref{alg:1} does not break down at iteration $k\leq \min\{n,p\}$.
Then the solution $y_k$ to \eqref{pro2} is
\begin{align}\label{pro3}
y_k=\arg\min_{y\in \mathbb{R}^{k}}||B_ky-\beta_1e_1||=\beta_1 B_k^{\dagger}e_1.
\end{align}	
\end{theorem}

{\em Proof.}
Let $\widetilde{y}=\bar{B}_ky$. Then
under the assumption on Algorithm~\ref{alg:1}, $\bar{B}_k$ is nonsingular.
Therefore, \eqref{pro2} is equivalent to
\begin{align*}
\min||\widetilde{y}|| \ \ \
\mbox{{\rm subject to}}\ \ \ \widetilde{y}\in \{\widetilde{y}|\
||(B_k\bar{B}_k^{-1})\widetilde{y}-\beta_1 e_1||=\min\}.
\end{align*}
Notice that $B_k$ is of column full rank, so is
$B_k\bar{B}_k^{-1}$. As a result, we have
\begin{align*}
\widetilde{y}_k &= \beta_1(B_k\bar{B}_k^{-1})^{\dagger}e_1\\
&=\beta_1(\bar{B}_kB_k^{\dagger})e_1
\end{align*}
with the second relation holding because
$B_k$ is of column full rank and $\bar{B}_k$ is nonsingular.
Then the solution $y_k$ to \eqref{pro2} is
\begin{align*}
y_k&= \bar{B}_k^{-1}\widetilde{y}_k
=\beta_1\bar{B}_k^{-1}\left((\bar{B}_kB_k^{\dagger})e_1\right)
=\beta_1B_k^{\dagger}e_1. \qquad \endproof
\end{align*}

\eqref{pro3} indicates that $y_k$ is simply the solution to the ordinary least
squares problem $\min_{y}||B_ky-\beta_1e_1||$ and $\bar{B}_k$ is
not invoked. Let
\begin{equation}
B_k=Q_kR_k
\end{equation}
be the compact QR factorization of $B_k$, which can be computed by exploiting
Givens rotations at cost of $\mathcal{O}(k)$ flops.
From \eqref{pro3}, we obtain
\begin{equation}
y_k=\beta_1R_k^{-1}Q_k^Te_1
\end{equation}
at cost of $\mathcal{O}(k)$ flops; see \cite{lsqr1982} for details.

Next, we consider the efficient computation of the
residual norm $\|Ax_k-b\|$ and the semi-norm $\|Lx_k\|$.

\begin{theorem}
Let the matrices $U_{k+1}$,
$\widehat{U}_k$, $B_k$ and $\bar{B}_k$ be defined in
\eqref{bid2} and \eqref{bid22}.
Then
\begin{align}\label{res}
\|Ax_k-b\|&=\|B_ky_k-\beta_1e_1\|,\\
\|Lx_k\|&=\|\bar{B}_ky_k\|.\label{sol}
\end{align}
\end{theorem}

{\em Proof.}
Notice $x_k=Z_ky_k$, and exploit \eqref{bid1} and \eqref{bid2}. We obtain
\begin{align*}
Ax_k=AZ_ky_k=U_{k+1}B_ky_k.
\end{align*}
Since $U_{k+1}$ is column orthonormal,
it is direct to derive \eqref{res} by the orthogonal invariance
of the 2-norm. Similarly, we have
\begin{align}\label{Lxx}
Lx_k=LZ_ky_k=\widehat{U}_k\bar{B}_ky_k.
\end{align}
Since $\widehat{U}_k$ is column orthognormal,
we have \eqref{sol}.
\qquad\endproof

This theorem shows that, by making use of structures of $B_k$
and $\bar{B}_k$, both $\|Ax_k-b\|$ and $\|Lx_k\|$
can be computed very efficiently at cost of $\mathcal{O}(k)$ flops
without forming $x_k$ explicitly. We only need to compute $x_k$ by
solving \eqref{xk2} when $x_k$ is accepted as the best regularized
solution.

\section{Regularization properties of the joint bidiagonalization
based algorithm}

Now we analyze our algorithm, establish some important results, and get
insight into its regularizing effects. Let $\widetilde{w}= Rx$. Then
by \eqref{qrfact}, we have
\begin{equation}\label{equiv}
\min_x\|Ax-b\|=\min_{\widetilde{w}} \|Q_A \widetilde{w}-b\|.
\end{equation}
First, it is direct to establish the following result, similar to Theorem 4.3 and
(3.4) in \cite{kilmer2007}.

\begin{lemma}\label{kry}
Let $x_k$ be the regularized solution obtained by our algorithm. Then
\begin{equation}\label{lsqrqa}
x_k=R^{-1} \widetilde{w}_k,\ \widetilde{w}_k=
\arg \min_{\widetilde{w}\in \mathcal{K}_k}\|Q_A \widetilde{w}-b\|,
\end{equation}
where $\mathcal{K}_k$ is the $k$ dimensional Krylov subspace
\begin{equation}\label{krylov}
\mathcal{K}_k={\rm span}\{Q_A^Tb,(Q_A^TQ_A)Q_A^Tb,\ldots,(Q_A^TQ_A)^{k-1}Q_A^Tb\},
\end{equation}
and the solution subspace
\begin{equation}\label{solspace}
{\rm span}\{Z_k\}=R^{-1}\mathcal{K}_k={\rm span}\{G(C^T C)^iC^TP_A^T b\}_{i=0}^{k-1},
\end{equation}
the $k$-dimensional Krylov subspace generated by the starting vector
$GC^TP_A^Tb$ and the matrix $C^TC$, where $C$ is defined by \eqref{csd}.
\end{lemma}

{\em Proof.} Write $\widetilde{w}=V_k y$, where $V_k$ is generated by
\eqref{diag1} and
$span\{V_k\}=\mathcal{K}_k$. Then from \eqref{diag1}
we obtain
$$
\min_{\widetilde{w}\in \mathcal{K}_k}\|Q_A \widetilde{w}-b\|=\min_y \|B_ky-\beta_1 e_1\|.
$$
Let $y_k=\arg \min_y \|B_ky-\beta_1 e_1\|$.
By the definition \eqref{vk} of $Z_k$, we have $x_k=Z_ky_k=R^{-1}V_ky_k=R^{-1}
\widetilde{w}_k$.  A direct justification
using \eqref{csd} and \eqref{gsvd} shows \eqref{solspace}.
\qquad\endproof

To present our main theoretical result and make an insightful
analysis on the regularizing effects of the proposed algorithm,
we need to make some necessary preparations and notation changes. Firstly, for
the regularization matrix $L\in \mathbb{R}^{p\times n}$ of rank $\min\{n,p\}$,
from the SVD \eqref{csd} of $Q_A$ and $Q_L$ and the labeling orders
\eqref{svda} and \eqref{svdl}, it is obvious that the singular values $c_i$ and
$s_i$ must satisfy
$$
0<c_i,\ s_i <1,\ i=1,2,\ldots,\min\{n,p\}.
$$
If $p\geq n$, we retain the notation \eqref{csd} and \eqref{svda}, and have
$1>c_1\geq c_2\geq \cdots \geq c_n>0$. If $p<n$,
different from \eqref{svda}, we relabel the $c_i$ and use the new notation
\begin{equation}\label{neworder}
1=c_1>c_2\geq c_3\geq\cdots\geq c_{p+1},
\end{equation}
where $c_1=1$ is the largest singular value of $Q_A$ with the multiplicity
$n-p$. That is, we reassign the indices $i$ of $c_i$ defined
by  \eqref{svda} to $i+1$, $i=1,2,\ldots,p$,
and shift the largest singular value of $Q_A$  to $c_1=1$ with multiplicity $n-p$
in \eqref{svda}. Correspondingly, we permute the columns of $P_A$ and $W$
in $Q_A$ defined by \eqref{csd}, and $G$ defined by \eqref{gsvd}
by moving their respective
last $n-p$ columns to the first ones and renaming
\begin{align*}
P_A:&=(P_{1,A},p_{2,A},\ldots,p_{m-n+p+1,A}),\\
W:&=(W_1,w_2,\ldots,w_{p+1}),\\
G:&=(G_1,g_2,\ldots,g_{p+1}).
\end{align*}
With the new notation, we have the range
$\mathcal{R}(G_1)=\mathcal{N}(L)$, i.e., $LG_1=0$.
Keep in mind that if $p\geq n$ then
$P_A$, $W$ and $G$ remain the same as in \eqref{csd} and \eqref{gsvd}, and
$\mathcal{N}(L)=\{\mathbf{0}\}$.

Secondly, it is well known \cite{bjorck96}
that the Lanczos bidiagonalization method for computing the singular values $c_i$
of $Q_A$ with the starting vector $b/\|b\|$ mathematically
amounts to the symmetric Lanczos method for computing the eigenvalues $c_i^2$ of
$Q_A^TQ_A$ with the starting vector $Q_A^Tb/\|Q_A^T b\|$. It is remarkable that
the symmetric Lanczos method works on $Q_A^TQ_A$ as if $Q_A^TQ_A$
has only simple eigenvalues $c_i^2$ \cite{parlett}. As a consequence,
the Lanczos bidiagonalization method works on $Q_A$ as if
the singular values $c_i$ of $Q_A$ are all simple. Notice that the singular values
$\widetilde{c}_i$, called the Ritz values, of the projected matrix $B_k$ are
always simple provided that the Lanczos bidiagonalization process does
not break down at step $k$. The Lanczos bidiagonalization method uses
the $\widetilde{c}_i$ as approximations to the $k$ {\em distinct} singular values
$c_i$ of $Q_A$. For a rigorous and complete derivation and many details,
we refer to \cite{jia18c}.

Next we establish an attractive and desired property that the regularized solution
$x_k$ has a filtered GSVD expansion and
is {\em explicitly} expressed in the generalized right singular vector
basis $\{g_i\}_{i=1}^n$  of $\{A,L\}$.

\begin{theorem}\label{filter}
Assume that the $c_i$ are labeled as \eqref{neworder} and simple for $p<n$,
the matrices $P_A,\ W$ and $G$ defined as above, and $g_{\perp}$ defined
by \eqref{second}. Then
\begin{equation}\label{xkfilter}
x_k=f_1^{(k)} g_{\perp}+
\sum_{i=2}^{p+1} f_i^{(k)}\frac{p_{i,A}^Tb}{c_i} g_i,
\ k=1,2,\ldots, n,
\end{equation}
where the filters
\begin{equation}\label{filt}
f_i^{(k)}=1-\prod_{j=1}^{k}\frac{\widetilde{c}_j^2-c_i^2}{\widetilde{c}_j^2}, \
i=1,2,\ldots, p+1.
\end{equation}
\end{theorem}

{\em Proof.}
Notice that the LSQR algorithm starting with $u_1=b/||b\|$ applied to \eqref{lsqrqa}
is mathematically equivalent to the conjugate gradient (CG) method applied to
the normal equation $Q_A^TQ_Aw=Q_A^Tb$ of \eqref{lsqrqa} with the starting vector
$w_0=0$.
Let $w_{ls}=Q_A^{\dagger}b$ be the solution to
$\min_{\widetilde{w}} \|Q_A \widetilde{w}-b\|$.
Then by the SVD \eqref{csd} of $Q_A$ and the
notation \eqref{neworder} we obtain
\begin{equation}\label{wls}
w_{ls}=W_1P_{1,A}^Tb+\sum_{i=2}^{p+1}\frac{p_{i,A}^Tb}{c_i}w_i,
\end{equation}
where the first term is the sum of the $n-p$ SVD
components of $Q_A$ corresponding to the largest singular value $c_1=1$
with the multiplicity $n-p$.

With our notation and \eqref{wls}, keep in mind
a well-known result (cf., e.g., \cite[Property 2.8]{vorst86}) on the CG iterates
that states
\begin{equation}\label{ritz}
\widetilde{w}_k=(I-q_k(Q_A^TQ_A))w_{ls},
\end{equation}
where $q_k(\mu)$ is the $k$-th residual polynomial of CG at iteration $k$
and $q_k(0)=1$,
whose roots are the Ritz values $\widetilde{c}_j^2$ of $Q_A^TQ_A$ with respect
to $\mathcal{K}_k$ defined by \eqref{krylov}, i.e.,
$$
q_k(c_i^2)=\prod_{j=1}^{k}\frac{\widetilde{c}_j^2-c_i^2}{\widetilde{c}_j^2},\
i=1,2,\ldots,p+1.
$$
Substituting \eqref{wls} into \eqref{ritz} yields
\begin{equation}\label{wkform}
\widetilde{w}_k=f_1^{(k)}W_1P_{1,A}^Tb+
\sum_{i=2}^{p+1} f_i^{(k)}\frac{p_{i,A}^Tb}{c_i} w_i,
\ k=1,2,\ldots, p+1
\end{equation}
with $f_i^{(k)}$ defined by \eqref{filt}.

Recall from \eqref{qrfact}--\eqref{gsvd} that $G=R^{-1}W$.
It then follows that $G_1=R^{-1}W_1$ and
$g_i=R^{-1}w_i,\ i=2,3,\ldots,p+1$. From \eqref{lsqrqa}, since
$x_k=R^{-1}\widetilde{w}_k$, premultiplying \eqref{wkform} by $R^{-1}$
establishes \eqref{xkfilter} by noticing that
$G_1P_{1,A}^Tb=g_{\perp}$ in our new notation.
\qquad\endproof

If $p\geq n$, then $g_{\perp}=0$ in \eqref{second},
the first term is zero in \eqref{xkfilter}, and the second term
becomes
$$
x_k=\sum_{i=1}^n f_i^{(k)}\frac{p_{i,A}^Tb}{c_i} g_i,
\ k=1,2,\ldots, n.
$$
In this case, \eqref{xkfilter} is a filtered GSVD expansion similar to
\eqref{tiksol}. If $p<n$, the first term
$$
f_1^{(k)} g_{\perp} \in \mathcal{N}(L)
$$
in \eqref{xkfilter},
which resembles the term
$g_{\perp}$ in \eqref{tiksol} and \eqref{tgsvdsol}. On the other hand,
$$
\sum_{i=2}^{p+1} f_i^{(k)}\frac{p_{i,A}^Tb}{c_i} g_i
$$
in \eqref{xkfilter} corresponds to the first term in \eqref{tiksol}
by noticing that in our notation the indices $i+1$ in the sum
correspond to the indices $i$ in \eqref{tiksol}.
A difference is that the general-form Tikhonov regularization solution
\eqref{tiksol} and TGSVD solution \eqref{tgsvdsol} do not affect $g_{\perp}$,
while our algorithm multiplies
it by a factor $f_1^{(k)}$. Nonetheless, $f_1^{(k)}\rightarrow 1$
as $\widetilde{c}_1$ converges to $c_1$;
since $g_{\perp}$ and $f_1^{(k)} g_{\perp}$ lie
in $\mathcal{N}(L)$, they
have no effect on $Lx_{\lambda}$, $Lx_k^{tgsvd}$ and $Lx_k$.

It is known from \cite[Theorem 2.1.1, p.23]{hansen98} that
the $c_i$ decay like the singular values $\sigma_i$ of $A$ when the matrix
$(A^T,L^T)^T$ is well conditioned, which is true provided that $L$ is well
conditioned, as is usually the case in practical applications.
In the meantime, notice from \eqref{csd} and \eqref{gsvd} that
$Q_A$ and $A$ share the same $P_A$ and the problems
$\min_{x}\|Ax-b\|$ and $\min_{w}\|Q_Aw-b\|$ have the same right-hand side $b$.
Therefore, the two problems satisfy the same discrete Picard condition.

Furthermore, as has been proved in \cite{jiarank},
since the $c_i$ decay and are clustered
at zero, the singular values of $B_k$
converge to the large singular values $c_i$ of $Q_A$ in natural order
for severely and moderately ill-posed problems until the occurrence of
semi-convergence of LSQR for solving $\min_{\widetilde{w}} \|Q_A \widetilde{w}-b\|$.
From \eqref{xkfilter} and \eqref{filt}, it is easily justified
that $f_i^{(k)}\approx 1$
for $i=1,2,\ldots,k$ and  $f_i^{(k)}\approx 0$ for $i=k+1,\ldots,p+1$
when the $k$ Ritz values $\widetilde{c}_j$ approximate the large
singular values of $Q_A$ in natural order; we refer the reader to
\cite[pp.~146-148]{hansen98} and \cite{jia16,jiarank}
for more details. This means that
$x_k$ mainly contains the first $k$ dominant GSVD components of $\{A,L\}$
and filters the others corresponding to the small generalized
singular values until the semi-convergence of the proposed joint
bidiagonalization based method. Precisely,
with the equivalence \eqref{equiv} and $x_k=R^{-1}\widetilde{w}_k$,
adapted the results of \cite{jiarank} to our current
context, this theorem shows that the proposed joint bidiagonalization
based method exhibit typical semi-convergence at some iteration $k^*$:
$x_k$ and $Lx_k$ converge to $x_{true}$ and $Lx_{true}$ for
$k\leq k^*$ and afterwards they are deteriorated by the noise $e$ and diverge
for $k>k^*$. Therefore,
the iteration number $k$ plays the role of the regularization parameter,
and the semi-convergence of the joint bidiagonalization based
algorithm occurs at iteration $k^*$, which is such that
$\|L(x_{k^*}-x_{true})\|$ is minimal over all $k=1,2,\ldots,\min\{n,p\}$.

\section{The determination of the optimal regularization parameter $k^*$}

For our joint bidiagonalization based algorithm,
since the residual norm $\|Ax_k-b\|=\|B_ky_k-\beta_1e_1\|$ monotonically
decreases and the semi-norm $\|Lx_k\|=\|\bar{B}_ky_k\|$
monotonically increases practically with respect to $k$,
the L-curve criterion and the discrepancy principle suit well for
a practical determination of $k^*$.  We plot the curve
$$\left(\log(\|B_ky_k-\beta_1e_1\|),\log(\|\bar{B}_ky_k\|)\right)$$
and then determine $k$ at its overall corner as an estimate of $k^*$.
This is routine, and we do not repeat the determination procedure;
see \cite{hansen98,hansen2007,hansen10}.

If $\|e\|$ or its accurate estimate is known in advance,
the discrepancy principle \cite{engl00,hansen98,hansen10} is the simplest and a
reliable choice. We stop the algorithm at the first iteration $k$
satisfying
\begin{equation}\label{disc}
\|Ax_k-b\|= \|B_ky_k-\beta_1e_1\|\leq \tau \|e\|
\end{equation}
with $\tau> 1$ slightly, e.g., $\tau=1.1$ or smaller.
We then use such $k$ as an estimate of the optimal regularization parameter $k^*$.
We mention that a $\tau>1$ considerably, e.g., $\tau=2$,
is generally unsafe and may underestimate $k^*$ substantially.

Embedded with the above parameter-choice methods, we can now present our joint
bidiagonalization based algorithm, called JBDQR and named Algorithm \ref{alg:2}.

\begin{algorithm}[htb]
\caption{(JBDQR) Given $A\in \mathbb{R}^{m\times n}$ and
$L\in \mathbb{R}^{p\times n}$, solve \eqref{tik2} and compute the regularized
solution $x_{k^*}$ at semi-convergence.}
\begin{algorithmic}[1]\label{alg:2}

\STATE Starting with $k=1$, run Algorithm \ref{alg:1}, and obtain the small
projected problem \eqref{pro3}.

\STATE Compute the minimum 2-norm solution $y_k$ to \eqref{pro3}.

\STATE Compute $\|Ax_k-b\|$ and $\|Lx_k\|$ by the formulas \eqref{res} and \eqref{sol}.

\STATE Determine the optimal regularization parameter $k^*$ by the L-curve criterion
or check if the discrepancy principle \eqref{disc} is satisfied. If $k^*$ is
not found, set $k=k+1$, and update Algorithm~\ref{alg:1}.
Then go to Step 2.

\STATE After $k^*$ is determined, form the regularized solution $x_{k^*}$
by solving \eqref{xk2}.
\end{algorithmic}
\end{algorithm}

\section{Numerical examples}

In this section, we report numerical experiments to demonstrate that
our JBDQR algorithm works well and the best regularized solutions obtained
by it are at least as accurate as those obtained by the hybrid
one proposed by Kilmer {\em et al.} \cite{kilmer2007}
and can be substantially more accurate than the latter ones.
We also compare the optimal regularization parameters determined
by the L-curve criterion and the discrepancy principle with the true
optimal $k^*$.

We choose some one dimensional examples from the
regularization toolbox \cite{hansen2007} and
some two dimensional problems from the Matlab Image Processing Toolbox and
\cite{berisha,nagy2004};
see Table \ref{tab1}, where the two dimensional image deblurring problems
{\sf rice} and {\sf mri} are from
the Matlab Image Processing Toolbox. We denote the relative noise level
$$
\varepsilon = \frac{\|e\|}{\|b_{true}\|}.
$$
For the noise-free problems $Ax_{true}=b_{true}$ in Table~\ref{tab1},
we add a white noise $e$ with zero mean and a prescribed noise level
$\varepsilon$ to $b_{true}$ and form the noisy $b=b_{true}+e$.
To simulate exact arithmetic, the complete reorthogonalization is used
in Algorithm~\ref{alg:1}.

\begin{table}[h]
    \centering
    \caption{The description of test problems.}
    \begin{tabular}{lll}
     \hline
     Problem        & Description                                & Ill-posedness \\
     \hline
     {\sf shaw}     & One-dimensional image restoration model  \cite{hansen2007}  & severe\\
     {\sf baart}  & First kind Fredholm integral equation \cite{hansen2007} & severe\\
     {\sf heat}     & Inverse heat equation                 \cite{hansen2007}     & moderate\\
     {\sf deriv2}   & Computation of second derivative      \cite{hansen2007}     & moderate\\
     {\sf rice}     &Two dimensional image deblurring        & unknown \\
     {\sf mri}      & Two dimensional image deblurring            & unknown \\
     {\sf AtmosphericBlur30} &Two dimensional image deblurring \cite{berisha,nagy2004}   & unknown\\
     {\sf GaussianBlur422}   &Two dimensional image deblurring \cite{berisha,nagy2004}   & unknown\\
     \hline
   \end{tabular}
   \label{tab1}
\end{table}

We abbreviate Algorithm~\ref{alg:2} as {\sf JBDQR},
the hybrid one in \cite{kilmer2007}
using the GCV and WGCV parameter-choice methods as
{\sf JBDGCV} and {\sf JBDWGCV}, respectively.
Let $x_k^{reg}$ denote the regularized solutions obtained by each of the
algorithms. We use the relative error
\begin{equation}\label{relerr}
\frac{\|L(x_k^{reg}-x_{true})\|}{\|Lx_{true}\|}
\end{equation}
to plot the convergence curve of each algorithm with respect to $k$.
In the tables to be presented,
we will list the smallest relative errors and iteration steps used by
{\sf JBDGCV} and {\sf JBDWGCV}
in parentheses, the optimal iteration steps
$k^*$ at which the semi-convergence of {\sf JBDQR} occurs
in the parentheses and the
estimated ones for $k^*$ determined by the L-curve
criterion and the discrepancy principle \eqref{disc} as well as
the corresponding relative errors
in the parentheses. We use the Matlab function {\sf lsqr.m}
to solve \eqref{innerleast}, \eqref{inner} and \eqref{xk2}
with the default stopping tolerance $tol=10^{-6}$.

All the computations are carried out in Matlab R2015b 64-bit on
Intel Core i3-2120 CPU 3.30GHz processor and 4 GB RAM with the machine precision
$\epsilon_{\rm mach}= 2.22\times10^{-16}$ under the Miscrosoft
Windows 7 64-bit system.

\subsection{One dimensional case}

The test problems {\sf shaw} and {\sf baart} are severely ill-posed, and {\sf heat}
and {\sf deriv2} are moderately ill-posed.
For each of them  we use the code of \cite{hansen2007} to generate $A$,
$x_{true}$ and $b_{true}$.
We mention that {\sf deriv2} has three
kinds of right-hand sides, distinguished by the parameter "$example=1,2,3$".
we only report the results on the parameter "$example=2$"
since we have obtained very
similar results on the problem with "$example=1,3$".
In the experiments,  for {\sf shaw} and {\sf baart}, we take $m=n=1024$, and
for ill-posed problems {\sf heat} and {\sf deriv2}, we take $m=n=3000$;
purely for test purposes, we choose
\begin{equation}\label{l1}
L=L_1 = \left(
        \begin{array}{ccccc}
          1 & -1 &  &  &  \\
           & 1 & -1 &  &  \\
           &  & \ddots & \ddots &  \\
             &  &  & 1  & -1\\
        \end{array}
      \right)\in \mathbb{R}^{(n-1)\times n},
\end{equation}
which is a scaled discrete approximation of
the first derivative operator in the one dimensional case. We comment that for the
scaled discrete approximation of the second derivative operator,
we have observed very similar phenomena. Hence we only report the results on $L=L_1$.

\begin{table}[h]
\centering
  \caption{The relative errors and estimates for
  the optimal regularization parameters $k^*$
  by the L-curve criterion for the test problems with
  $L=L_1$.}
  \label{tab2}
  \centering{$\varepsilon=10^{-2}$}
  \begin{minipage}[t]{1\textwidth}
  \begin{tabular*}{\linewidth}{lp{2.4cm}p{2.4cm}p{2.4cm}p{2.4cm}p{2.4cm}p{2.4cm}}
  \toprule[0.6pt]
     	            &{\sf JBDGCV}    &{\sf JBDWGCV}    &{\sf JBDQR ($k^*$)}   &
     estimates for $k^*$\\ \midrule[0.6pt]
     {\sf shaw}     &0.5398(18)     &0.5398(18)      &0.2094(3)    &2(0.2126)	         \\
     {\sf baart}    &0.5574(7)      &0.5574(7)       &0.5405(2)    &3(0.5625)     \\
     {\sf heat}     &0.3758(60)     &0.3758(60)      &0.2186(13)   &5(0.3284)  \\
     {\sf deriv2}   &0.4270(60)     &0.4270(60)      &0.3363(4)    &2(0.3853) \\
    \bottomrule[0.6pt]
     \end{tabular*}\\[2pt]
      \end{minipage}
      \centering{$\varepsilon=10^{-3}$}
  \begin{minipage}[t]{1\textwidth}
  \begin{tabular*}{\linewidth}{lp{2.4cm}p{2.4cm}p{2.4cm}p{2.4cm}p{2.4cm}p{2.4cm}}
  \toprule[0.6pt]
     	            &{\sf JBDGCV}    &{\sf JBDWGCV}    &{\sf JBDQR ($k^*$)}   &
     estimates for $k^*$\\ \midrule[0.6pt]
     {\sf shaw}     &0.1930(13)     &0.1930(13)      &0.1732(5)    &2(0.1918)	         \\
     {\sf baart}    &0.5442(9)      &0.5442(9)       &0.5038(4)    &2(0.5376)     \\
     {\sf heat}     &0.1794(100)    &0.1794(100)     &0.1456(25)   &23(0.1485)  \\
     {\sf deriv2}   &0.3884(60)     &0.3884(60)      &0.2635(10)   &8(0.3161) \\
    \bottomrule[0.6pt]
    \end{tabular*}\\[2pt]
    \end{minipage}
  \centering{$\varepsilon=10^{-4}$}
  \begin{minipage}[t]{1\textwidth}
  \begin{tabular*}{\linewidth}{lp{2.4cm}p{2.4cm}p{2.4cm}p{2.4cm}p{2.4cm}p{2.4cm}}
  \toprule[0.6pt]
     	           &{\sf JBDGCV}    &{\sf JBDWGCV}    &{\sf JBDQR ($k^*$)}    &estimates for $k^*$\\ \midrule[0.6pt]
     {\sf shaw}    &0.1664(14)   &0.1664(14)    &0.1378(8)     &8(0.1378)	\\
     {\sf baart}   &0.5346(9) 	 &0.5346(9)	    &0.4136(5)	   &3(0.5354)     \\
     {\sf heat}    &0.1360(100)	 &0.1360(100)	&0.1275(37)	   &35(0.1283)   \\
     {\sf deriv2}  &0.2916(60)	 &0.2916(60)	&0.2452(15)	   &12(0.2606)  \\
    \bottomrule[0.6pt]
    \end{tabular*}\\[2pt]
  \end{minipage}
\end{table}

\begin{table}[h]
\centering
  \caption{The relative errors and estimates for the optimal regularization
  parameters $k^*$
  by the discrepancy principle for the test problems with
  $L=L_1$.}
  \label{addtab1}
  \centering{$\varepsilon=10^{-2}$}
  \begin{minipage}[t]{1\textwidth}
  \begin{tabular*}{\linewidth}{lp{2.4cm}p{2.4cm}p{2.4cm}p{2.4cm}p{2.4cm}p{2.4cm}}
  \toprule[0.6pt]
     	            &$\tau=1.005$    &$\tau=1.1$    &$\tau=1.2$   &$\tau=2.0$\\ \midrule[0.6pt]
     {\sf shaw}     &0.3031(1)     &0.3031(1)      &0.3031(1)    &0.3031(1)	         \\
     {\sf baart}    &0.5421(1)     &0.5421(1)      &0.5421(1)    &0.5421(1)     \\
     {\sf heat}     &0.3152(6)     &0.3757(3)      &0.4629(2)    &0.5410(1)  \\
     {\sf deriv2}   &0.3853(2)     &0.4187(1)      &0.4187(1)    &0.4187(1) \\
    \bottomrule[0.6pt]
     \end{tabular*}\\[2pt]
      \end{minipage}
      \centering{$\varepsilon=10^{-3}$}
  \begin{minipage}[t]{1\textwidth}
  \begin{tabular*}{\linewidth}{lp{2.4cm}p{2.4cm}p{2.4cm}p{2.4cm}p{2.4cm}p{2.4cm}}
  \toprule[0.6pt]
     	            &$\tau=1.005$    &$\tau=1.1$    &$\tau=1.2$   &$\tau=2.0$\\ \midrule[0.6pt]
     {\sf shaw}     &0.1888(2)      &0.1888(2)       &0.1888(2)    &0.2338(1)	         \\
     {\sf baart}    &0.5376(2)      &0.5422(1)       &0.5422(1)    &0.5422(1)     \\
     {\sf heat}     &0.1669(20)     &0.2196(11)      &0.2377(10)   &0.3258(5)  \\
     {\sf deriv2}   &0.3398(6)      &0.4291(2)       &0.4291(2)    &0.4651(1) \\
    \bottomrule[0.6pt]
    \end{tabular*}\\[2pt]
    \end{minipage}
  \centering{$\varepsilon=10^{-4}$}
  \begin{minipage}[t]{1\textwidth}
  \begin{tabular*}{\linewidth}{lp{2.4cm}p{2.4cm}p{2.4cm}p{2.4cm}p{2.4cm}p{2.4cm}}
  \toprule[0.6pt]
     	          &$\tau=1.005$    &$\tau=1.1$    &$\tau=1.2$   &$\tau=2.0$\\ \midrule[0.6pt]
     {\sf shaw}    &0.1632(5)    &0.1882(3)    &0.1882(3)      &0.1906(2)	\\
     {\sf baart}   &0.5354(3) 	 &0.5354(3)	   &0.5400(2)	   &0.5437(1)     \\
     {\sf heat}    &0.1356(28)	 &0.1443(25)   &0.1473(24)	   &0.1745(19)   \\
     {\sf deriv2}  &0.2606(12)	 &0.3019(9)	   &0.3400(7)	   &0.3813(4)  \\
    \bottomrule[0.6pt]
    \end{tabular*}\\[2pt]
  \end{minipage}
\end{table}

In Table \ref{tab2}, we display the relative errors of the best regularized
solutions
by {\sf JBDQR}, {\sf JBDGCV} and {\sf JBDWGCV} with $L=L_1$ and
$\varepsilon=10^{-2},\ 10^{-3}, \ 10^{-4}$, respectively. As we can see
from the table, the best regularized
solutions by {\sf JBDQR} are at least as accurate as and can be
considerably more accurate than those by
{\sf JBDGCV} and {\sf JBDWGCV} for all the test problems; see, e.g.,
{\sf shaw}, {\sf heat} and {\sf deriv2} for $\varepsilon=10^{-2}$, and
{\sf deriv2} for $\varepsilon=10^{-3}$.
We observe from the table that for
each test problem the best regularized solution by {\sf JBDQR}
is correspondingly more accurate
and requires a bigger $k^*$ for a smaller $\varepsilon$.
All these are expected and justify that the smaller $\varepsilon$ is,
the better regularized solution is extracted, that is,
the more GSVD dominant components of $\{A,L\}$ are needed to form it.
Finally, for {\sf JBDQR}, we see that
for each problem and given $\varepsilon$, almost all
the regularization parameters $k^*$ determined by the
L-curve criterion are quite reliable and close to the true $k^*$
except {\sf shaw} for $\varepsilon=10^{-3}$. But we also find that
the the L-curve criterion
underestimates the true $k^*$ more or less, that is, the estimates for $k^*$ by
the L-curve criterion oversmooths the regularized solutions.

\begin{figure}
\begin{minipage}{0.48\linewidth}
  \centerline{\includegraphics[width=6.0cm,height=4cm]{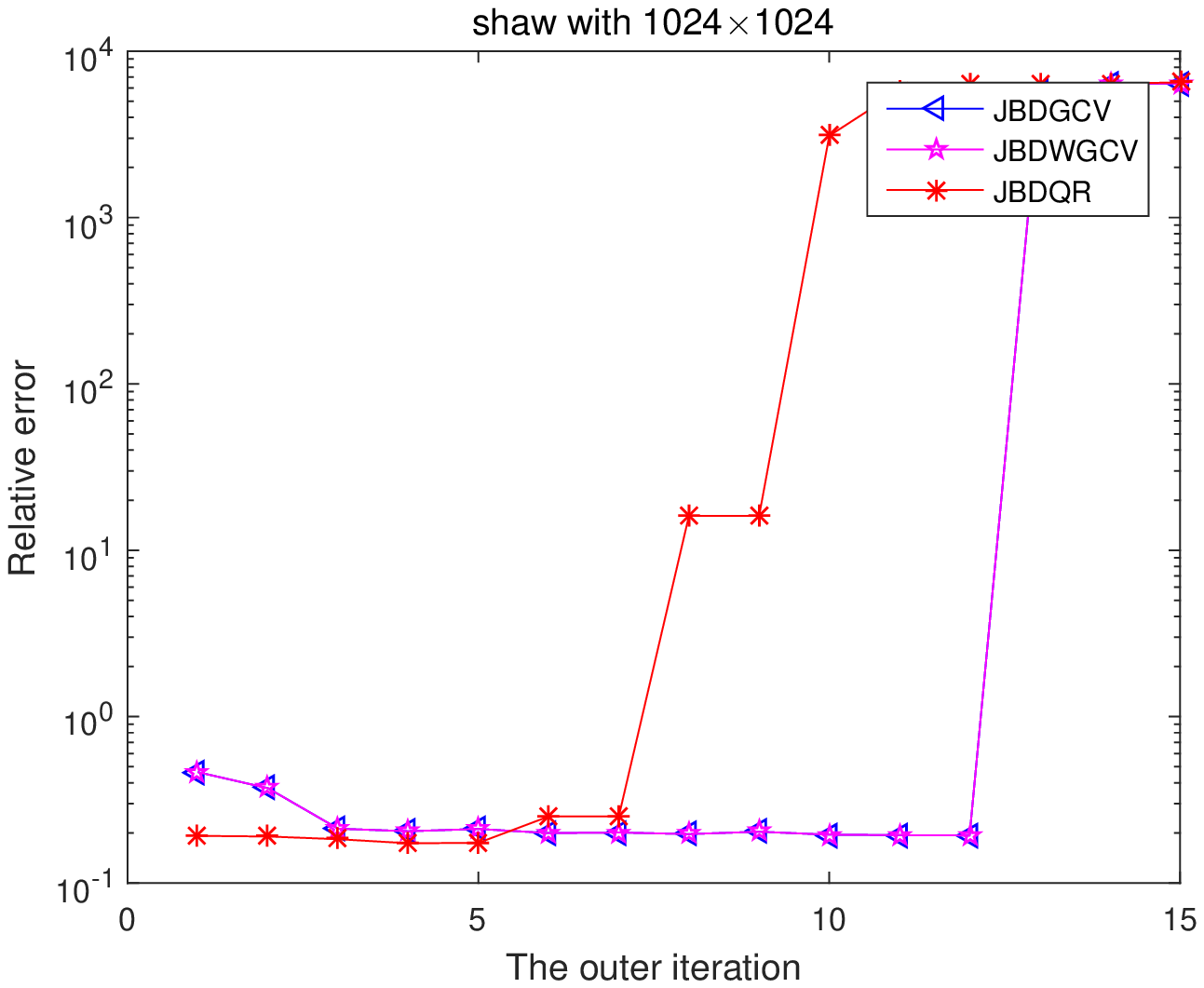}}
  \centerline{(a)}
\end{minipage}
\hfill
\begin{minipage}{0.48\linewidth}
  \centerline{\includegraphics[width=6.0cm,height=4cm]{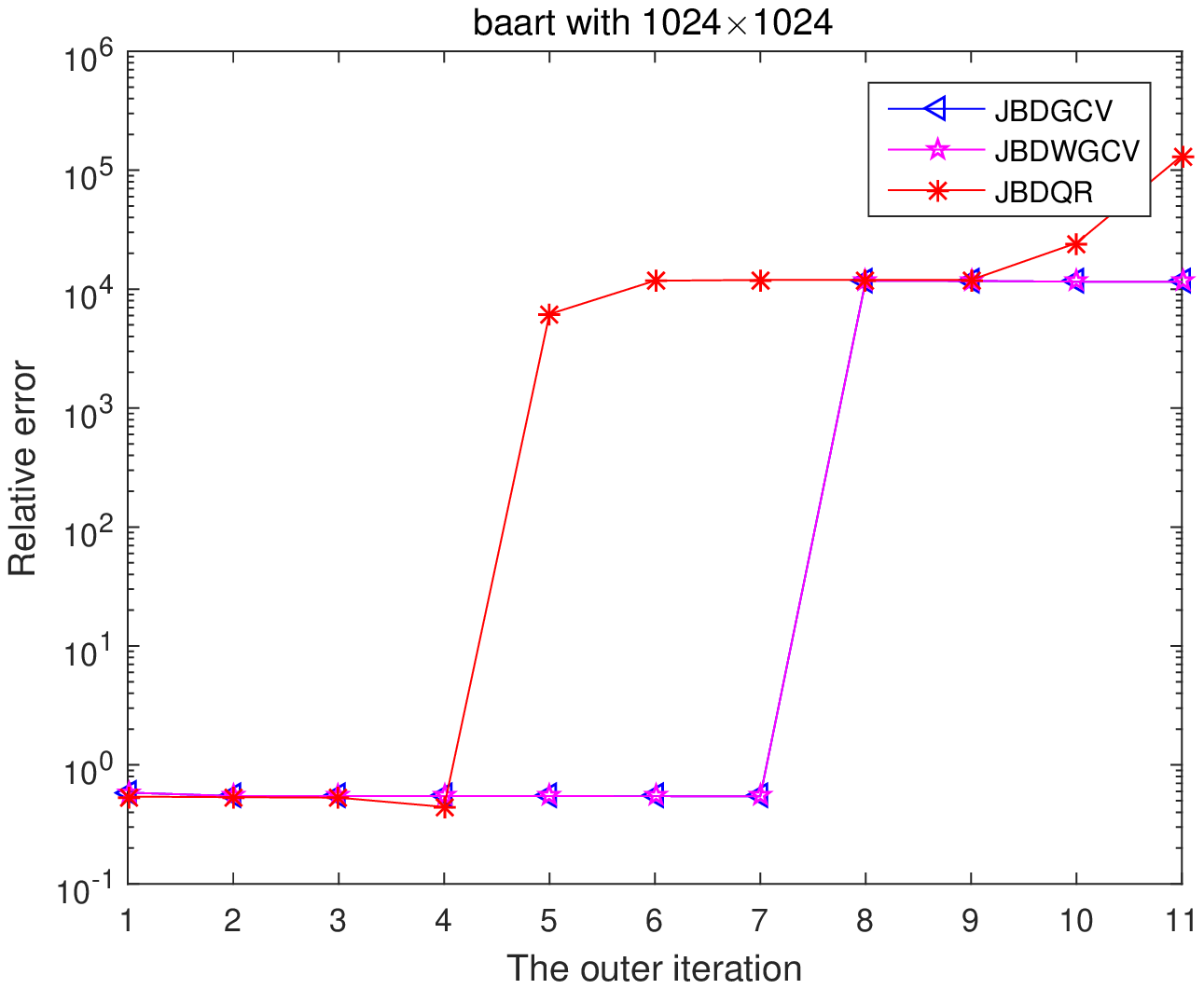}}
  \centerline{(b)}
\end{minipage}

\vfill
\begin{minipage}{0.48\linewidth}
  \centerline{\includegraphics[width=6.0cm,height=4cm]{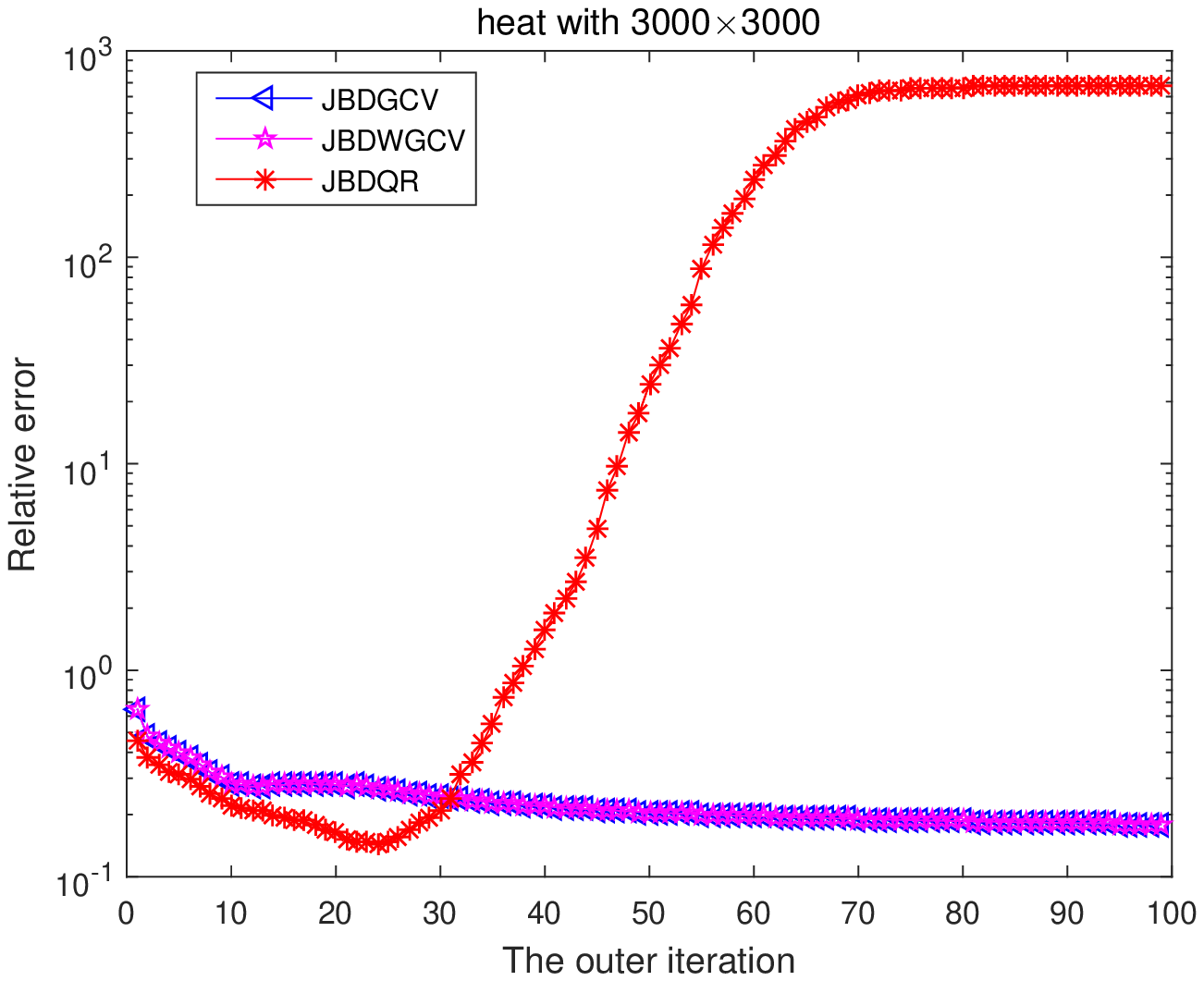}}
  \centerline{(c)}
\end{minipage}
\hfill
\begin{minipage}{0.48\linewidth}
  \centerline{\includegraphics[width=6.0cm,height=4cm]{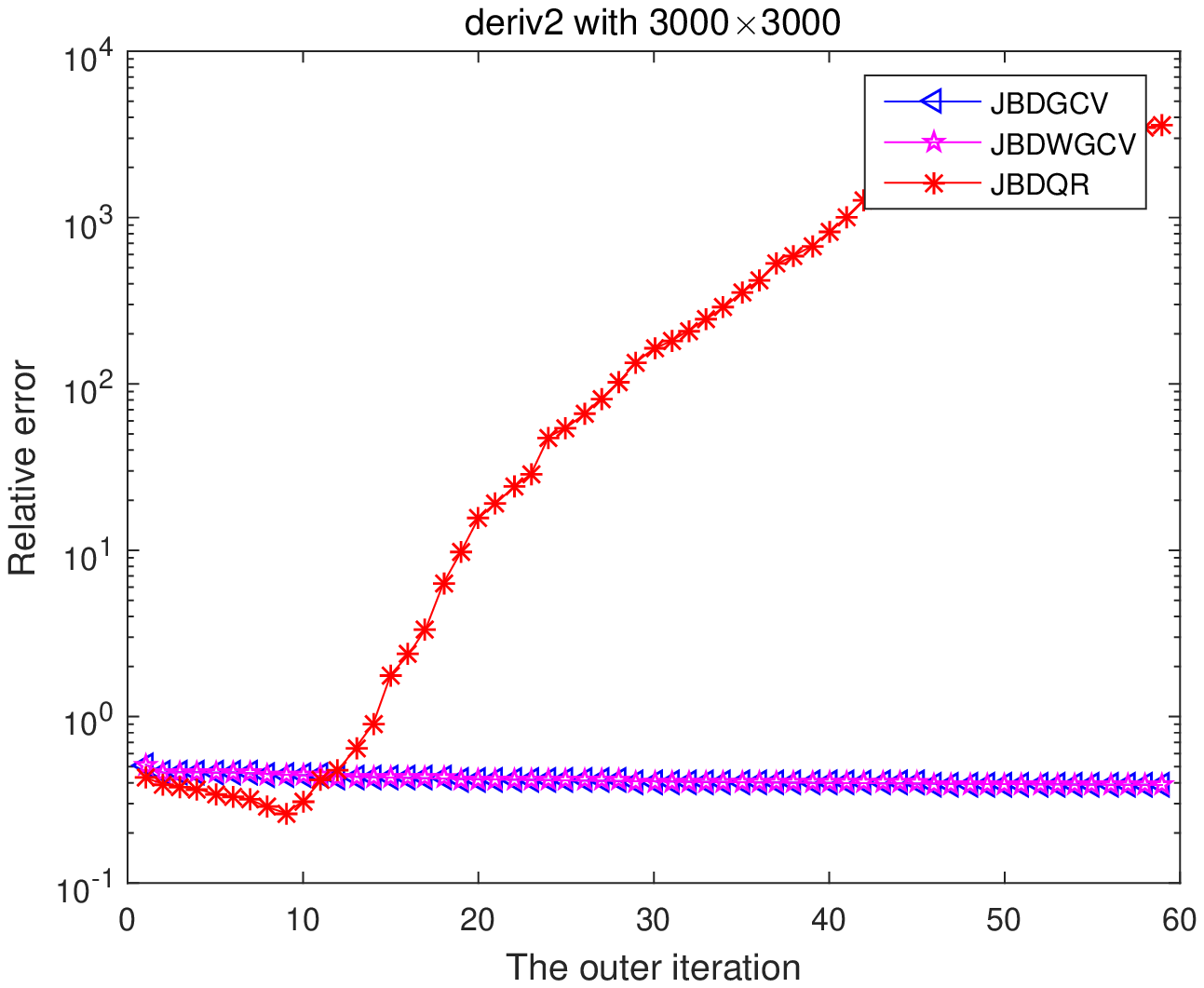}}
  \centerline{(d)}
\end{minipage}
\caption{The relative error of {\sf JBDQR}, {\sf JBDGCV} and {\sf JBDWGCV}
 with $L=L_1$ and $\varepsilon=10^{-3}$: (a) {\sf shaw};
 (b) {\sf baart}; (c) {\sf heat}; (d) {\sf deriv2}.}
\label{fig1}
\end{figure}

Figure \ref{fig1} depicts the convergence processes of
{\sf JBDQR}, {\sf JBDGCV} and {\sf JBDWGCV}
for $L=L_1$ and $\varepsilon=10^{-3}$.
We observe from the figure and Table~\ref{tab2}
that, in the most cases, the best regularized
solutions by {\sf JBDQR} are more accurate and can be considerably
more accurate than
those by {\sf JBDGCV} and {\sf JBDWGCV}. In addition, for
the severely ill-posed {\sf shaw} and {\sf baart} we find
that {\sf JBDGCV} and {\sf JBDWGCV} behave very similarly and the
convergence processes are almost indistinguishable. Remarkably, we see
that the regularized solutions obtained by them converge first,
then stabilize for a while, and finally diverge dramatically,
while, for {\sf heat} and {\sf deriv2}, they start to stabilize
after $k$ becomes large. We have also observed that the smaller
$\varepsilon$ is, the later they start to stabilize, though we do not
draw all the corresponding figures.
The phenomena for {\sf shaw} and {\sf baart} do not comply with the
expectation that the regularized solutions ultimately stabilize
as the subspace is expanded sufficiently large. The reason is
due to the fact that the discrete Picard conditions for the projected
problems are satisfied poorly as $k$ increases, as we have
argued in the introduction. In contrast, {\sf JBDQR} has always exhibited
the typical semi-convergence for all the problems, which justifies our theory.

\begin{figure}
\begin{minipage}{0.48\linewidth}
  \centerline{\includegraphics[width=6.0cm,height=4cm]{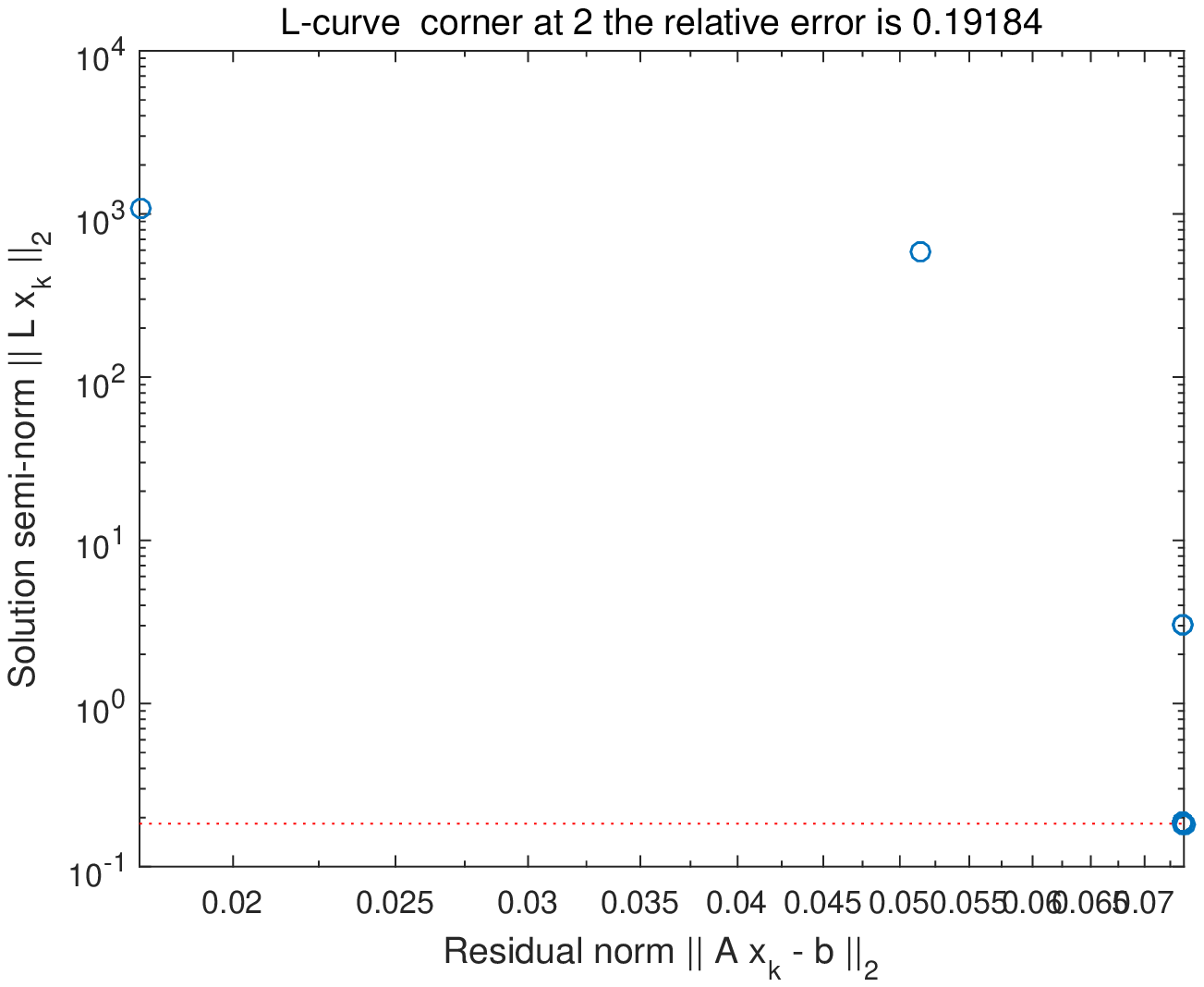}}
  \centerline{(a)}
\end{minipage}
\hfill
\begin{minipage}{0.48\linewidth}
  \centerline{\includegraphics[width=6.0cm,height=4cm]{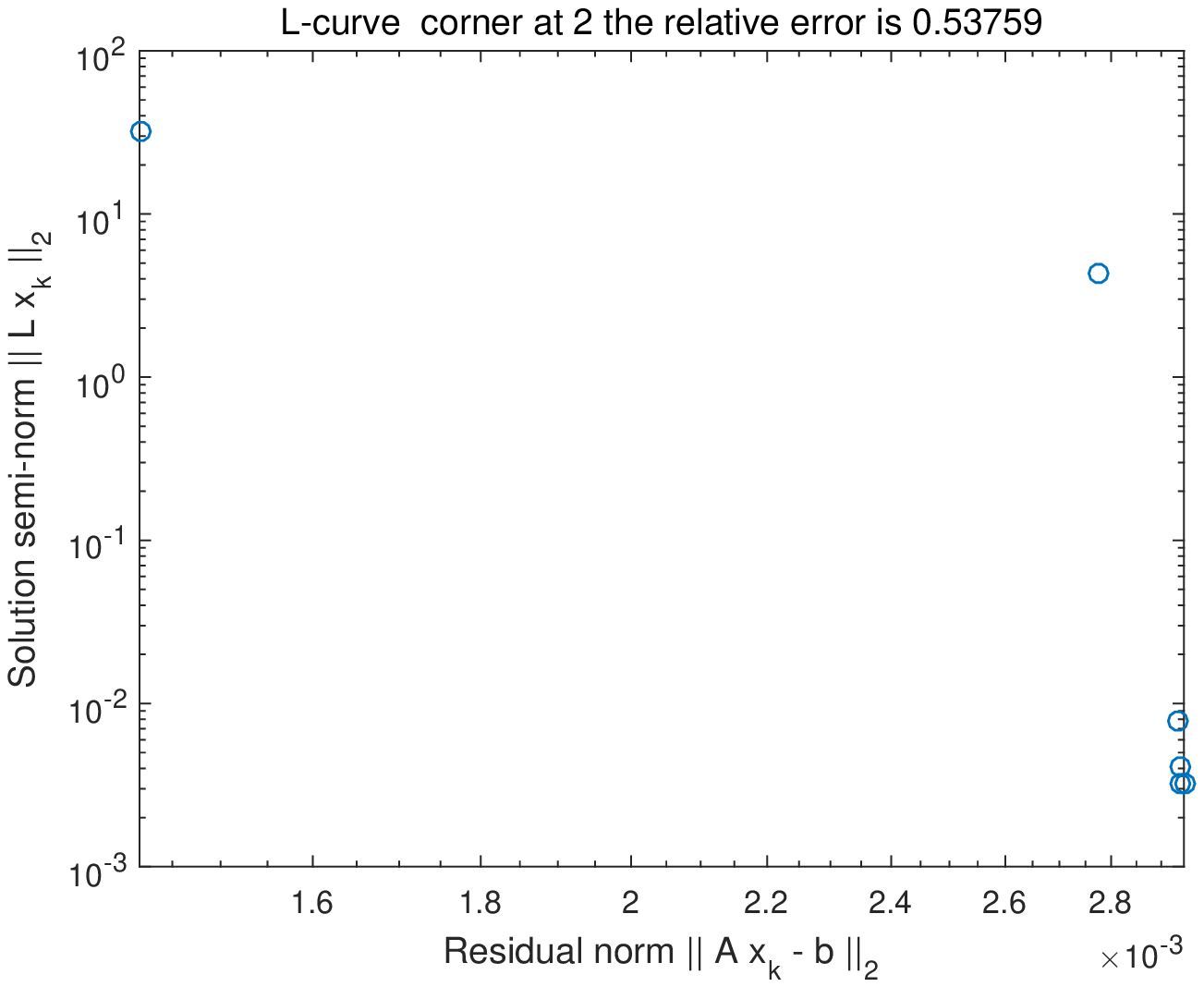}}
  \centerline{(b)}
\end{minipage}
\vfill
\begin{minipage}{0.48\linewidth}
  \centerline{\includegraphics[width=6.0cm,height=4cm]{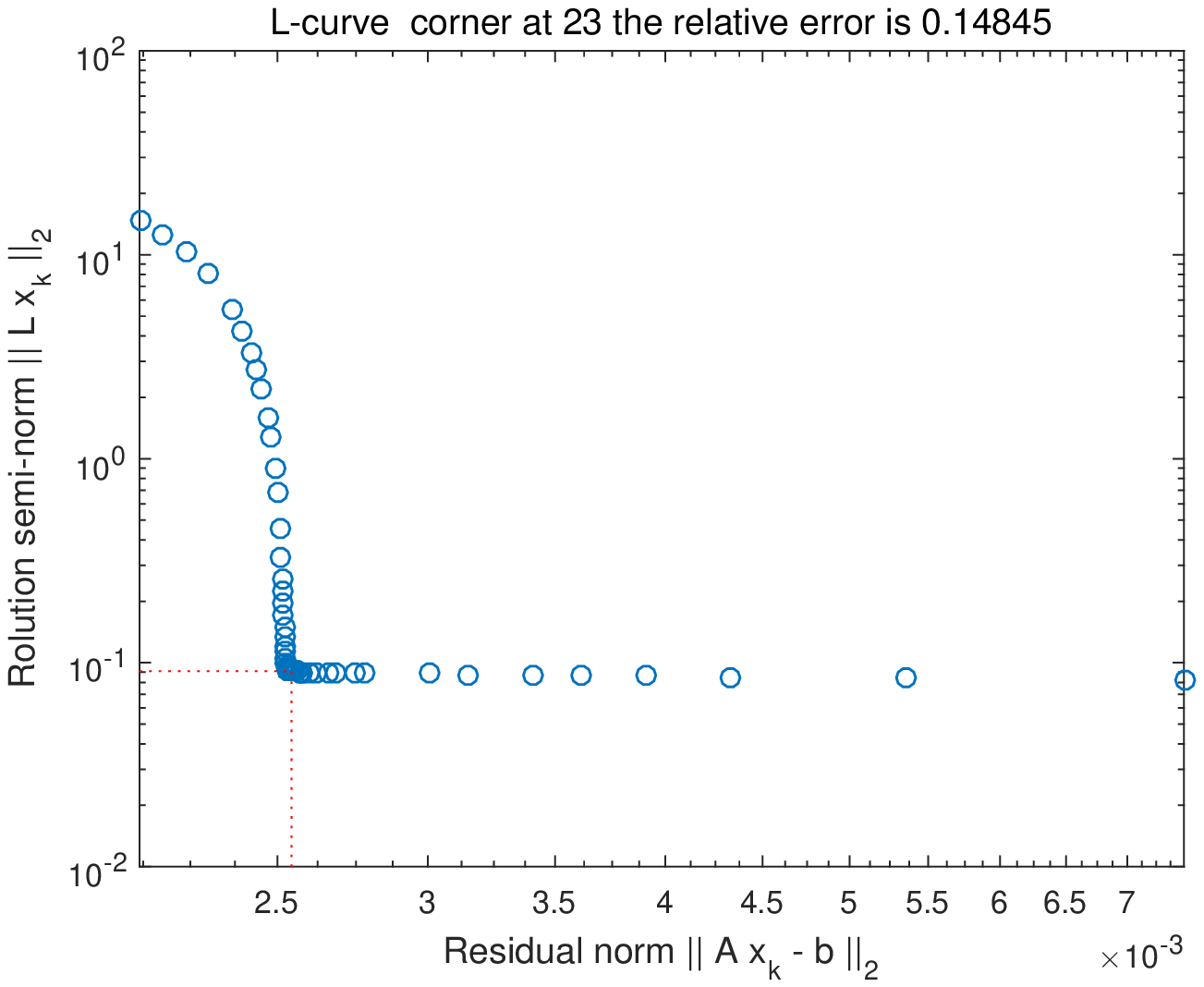}}
  \centerline{(c)}
\end{minipage}
\hfill
\begin{minipage}{0.48\linewidth}
  \centerline{\includegraphics[width=6.0cm,height=4cm]{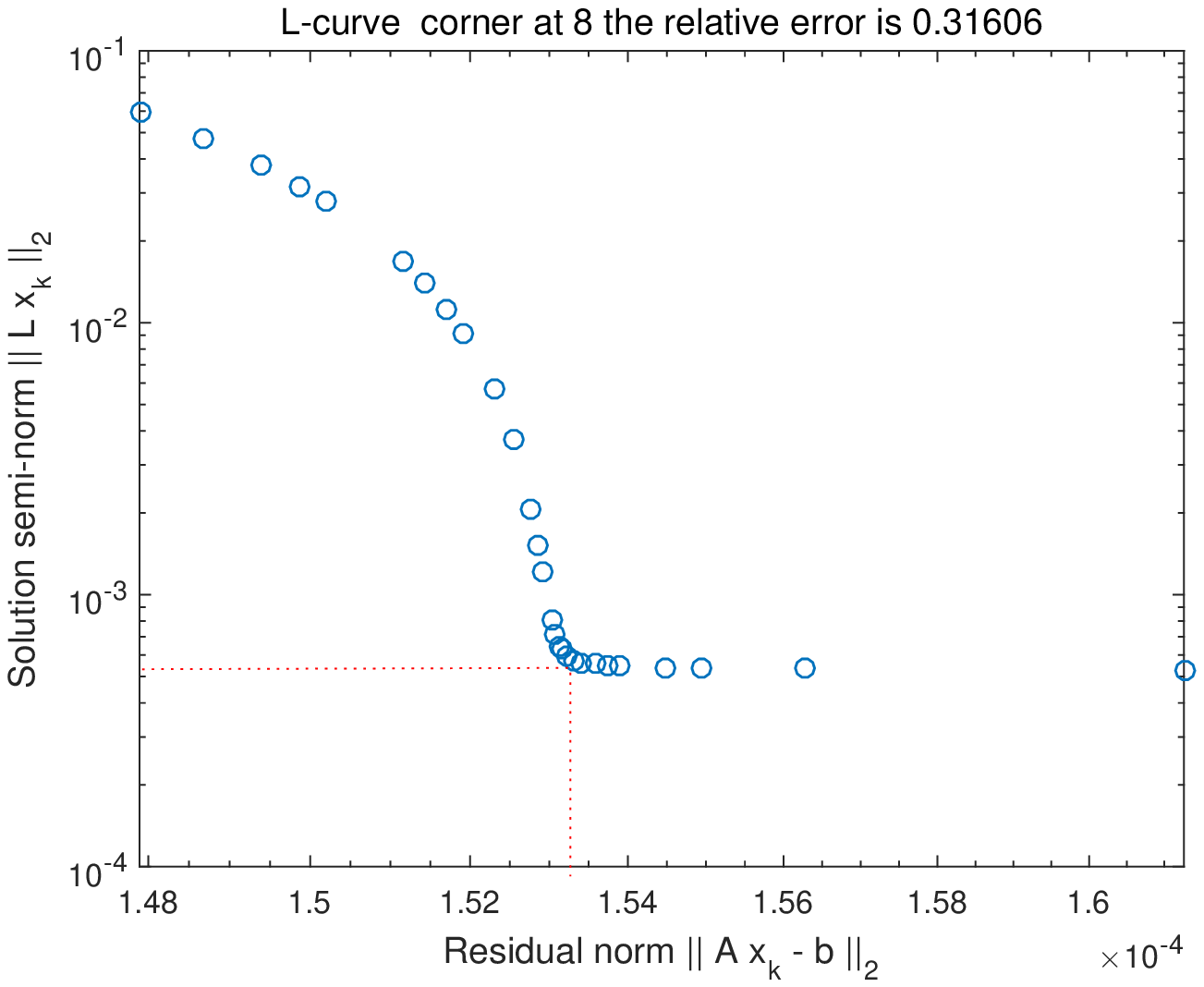}}
  \centerline{(d)}
\end{minipage}
\caption{The determination of $k^*$ by the L-curve criterion
by {\sf JBDQR} with $L=L_1$,
   $\varepsilon=10^{-3}$:
   (a) {\sf shaw};(b) {\sf baart}; (c) {\sf heat}; (d) {\sf deriv2}.}
\label{fig2}
\end{figure}

Figure \ref{fig2} depicts the L-curves given by {\sf JBDQR} with
$L=L_1$ and $\varepsilon=10^{-3}$. We use the function
{\sf l\_corner} in \cite{hansen2007} to determine the overall corner and
give an estimate $k^*$. We see that for the moderately
ill-posed problems {\sf heat} and {\sf deriv2} there
are much better "L" shape curves, which enable us to determine
the optimal $k^*$  more reliably and accurately
than those for the severely ill-posed problems {\sf shaw} and {\sf baart}. This
is because {\sf JBDQR} converges very fast and uses very few iterations to
achieve the semi-convergence for {\sf shaw} and {\sf baart}. Indeed,
the L-curve criterion does not work well for {\sf shaw}  and {\sf baart} with
$\varepsilon=10^{-3}$.

Since $\|e\|$ is known for the above test problems, we can
use the discrepancy principle \eqref{disc}
to estimate the optimal $k^*$.  Table~\ref{addtab1}
reports the results obtained, in which we have taken the
four $\tau=1.005, \
1.1, \ 1.2$ and $2.0$. Compared with the $k^*$ in Table~\ref{tab2},
we have found that the discrepancy principle always underestimate
$k^*$ and the problems are over-regularized.
We have observed that the reliable determination of $k^*$
critically depend on $\tau$, and the closer $\tau$ is to one, the more
reliable the estimates are. Particularly, except for {\sf shaw} and
{\sf baart} with $\varepsilon=10^{-2}$, the choice
$\tau=2$ is obviously very bad, and it produces very poor estimates for $k^*$
and leads to considerably less accurate regularized solutions than
$\tau=1.005$ does.

\subsection{Two dimensional case}

In this section, we test some two dimensional image deblurring problems
listed in Table~\ref{tab1}. The goal is to restore an image $x_{true}$
from a blurred and noisy image $b=b_{true}+e$.

We consider the problems {\sf rice} and {\sf mri} from
the Matlab Image Processing Toolbox.
The exact image $x_{true}$ of {\sf rice} is an $N\times N$
subimage and that of
{\sf mri} is the $15$th slice of the three dimensional
MRI image dataset which has $N\times N$ pixels.
The blurred operator $A$ is
a symmetric doubly Toeplitz PSF matrix and is of Kroneck product form
$A=(2\pi\sigma^2)^{-1}T\otimes T\in \mathbb{R}^{N^2\times N^2}$,
where $T\in\mathbb{R}^{N\times N}$ is a symmetric banded
Toeplitz matrix with half-bandwidth {\sf band} and
$\sigma$ controls the width of Gaussian PSF.
In what follows, we use
${\sf band}=16$, $\sigma=2$ and $N=128$. The size of {\sf rice}
and {\sf mri} is $m=n=128^2=16,284$.

We also consider the problems
{\sf AtmosphericBlur30} and {\sf GaussianBlur422} of $m=n=256^2=65,536$
from \cite{nagy2004}.
The blurring of {\sf AtmosphericBlur30} is caused by atmospheric turbulence,
and {\sf GaussianBlur422} is spatially invariant Gaussian blur.
The exact images are generated by the input command
``{\sf load  AtmosphericBlur30}''
and ``{\sf load  GaussianBlur422}'',
and the blurring operators are generated by the codes
{\sf psfMatrix(PSF,center, 'zero')} and
{\sf psfMatrix(PSF)} from \cite{nagy2004}, respectively.
We abbreviate {\sf AtmosphericBlur30} and {\sf GaussianBlur422} as
{\sf blur30} and {\sf blur422}, respectively.

For the experimental purpose, we choose the regularization matrix
\begin{equation}\label{l3}
L=\left(
    \begin{array}{c}
      I_N\otimes L_1 \\
      L_1 \otimes I_N \\
    \end{array}
  \right)\in\mathbb{R}^{N(N-1)\times N^2}
\end{equation}
with $L_1$ defined in \eqref{l1} and
$I_N$ the identity matrix of order $N$,
which is the scaled discrete approximation of the first derivative operator
in the two dimensional case incorporating
no assumptions on boundary conditions; see \cite[Chapter 8.1]{hansen10}.
The white noise
$e$ with zero mean are generated so that the relative noise levels $\varepsilon=
5\cdot 10^{-2}, \ 10^{-2}$ and $10^{-3}$, respectively.

\begin{table}[h]
  \centering
  \caption{The relative errors and estimates for $k^*$ by
  the L-curve criterion.}
  \label{tab3}
  \centering{$\varepsilon=5\cdot10^{-2}$}
  \begin{minipage}[t]{1\textwidth}
  \begin{tabular*}{\linewidth}{lp{2.0cm}p{2.0cm}p{2.2cm}p{2.4cm}p{2.4cm}p{2.5cm}}
  \toprule[0.6pt]
     	  &{\sf JBDGCV}    &{\sf JBDWGCV}    &{\sf JBDQR ($k^*$)} &estimates for
     $k^*$\\ \midrule[0.6pt]
     {\sf rice}           &0.8778(3)     &0.8736(3)     &0.8397(5)     &4(0.8411)\\
     {\sf rice(no $L$)}   &0.1175(3)     &0.1142(3)     &0.0950(4)      \\
     {\sf mri}            &0.9602(3)     &0.9498(4)     &0.8848(13)     &6(0.9007)\\
     {\sf mri(no $L$)}    &0.3066(3)     &0.2932(3)     &0.2324(13)      \\
     {\sf blur30}         &0.9835(3)     &0.9827(3)     &0.9124(16)    &5(0.9521)\\
     {\sf blur30(no $L$)} &0.5056(3)     &0.4999(3)     &0.3036(16)\\
     {\sf blur422}        &0.9459(9)     &0.9443(10)    &0.9109(62)    &24(0.9203)\\
     {\sf blur422(no $L$)}&0.2843(9)     &0.2823(10)    &0.2522(58)      \\
    \bottomrule[0.6pt]
    \end{tabular*}\\[2pt]
    \end{minipage}
     \centering{$\varepsilon=10^{-2}$}
      \begin{minipage}[t]{1\textwidth}
       \begin{tabular*}{\linewidth}{lp{2.0cm}p{2.0cm}p{2.4cm}p{2.4cm}p{2.4cm}p{2.4cm}}
        \toprule[0.6pt]
     	                 &{\sf JBDGCV}    &{\sf JBDWGCV}    &{\sf JBDQR ($k^*$)}  	&estimates
     for $k^*$\\ \midrule[0.6pt]
     {\sf rice}           &0.8372(7)     &0.8363(7)     &0.7774(23)    &11(0.7951)\\
     {\sf rice(no $L$)}   &0.0931(7)     &0.0927(7)     &0.0764(22)      \\
     {\sf mri}            &0.8923(13)    &0.8782(19)    &0.8421(51)    &25(0.8514)\\
     {\sf mri(no $L$)}    &0.2367(13)    &0.2258(19)    &0.2024(50)      \\
     {\sf blur30}         &0.9697(6)     &0.9603(9)     &0.7975(65)    &40(0.8243)\\
     {\sf blur30(no $L$)} &0.4238(5)     &0.3900(9)     &0.2095(65)      \\
     {\sf blur422}        &0.9459(9)     &0.9443(10)    &0.9109(62)    &24(0.9203)\\
     {\sf blur422(no $L$)}&0.2843(9)     &0.2823(10)    &0.2522(58)      \\
          \bottomrule[0.6pt]
          \end{tabular*}\\[2pt]
          \end{minipage}
          \centering{$\varepsilon=10^{-3}$}
      \begin{minipage}[t]{1\textwidth}
       \begin{tabular*}{\linewidth}{lp{2.0cm}p{2.0cm}p{2.4cm}p{2.4cm}p{2.4cm}p{2.4cm}}
        \toprule[0.6pt]
     	         &{\sf JBDGCV}  &{\sf JBDWGCV}  &{\sf JBDQR ($k^*$)}  &estimates for $k^*$\\ \midrule[0.6pt]
     {\sf rice}            &0.7638(38)     &0.7539(52)     &0.7136(166)      &145(0.7140)\\
     {\sf rice(no $L$)}    &0.0726(38)     &0.0705(52)     &0.0626(163)      \\
     {\sf mri}             &0.8305(101)    &0.8225(151)    &0.7949(451)     &293(0.7988)\\
     {\sf mri(no $L$)}     &0.1957(101)    &0.1917(151)    &0.1799(447)      \\
     {\sf blur30}          &0.9628(9)    &0.7984(75)     &0.5670(433)     &577(0.5907)\\
     {\sf blur30(no $L$)}  &0.3988(9)    &0.2060(75)     &0.1110(438)\\
     {\sf blur422}         &0.9137(59)     &0.9046(103)     &0.8736(542)    &284(0.8794)\\
     {\sf blur422(no $L$)} &0.2536(59)     &0.2471(103)     &0.2285(549)      \\
          \bottomrule[0.6pt]
          \end{tabular*}\\[2pt]
          \end{minipage}
\end{table}

\begin{table}[h]
  \centering
  \caption{The relative errors and estimates for $k^*$ by
  the discrepancy principle.}
  \label{addtab2}
  \centering{$\varepsilon=5\cdot10^{-2}$}
  \begin{minipage}[t]{1\textwidth}
  \begin{tabular*}{\linewidth}{lp{2.0cm}p{2.0cm}p{2.0cm}p{2.0cm}p{2.0cm}p{2.0cm}}
  \toprule[0.6pt]
     	               &$\tau=1.005$    &$\tau=1.1$    &$\tau=1.2$   &$\tau=2.0$\\ \midrule[0.6pt]
     {\sf rice}     &0.8462(3)     &0.8556(2)     &0.8556(2)    &0.8791(1)\\
     {\sf mri}      &0.9007(6)     &0.9156(4)     &0.9267(3)    &0.9653(1)    \\
     {\sf blur30}   &0.9181(12)    &0.9575(5)     &0.9700(3)    &0.9849(1)\\
     {\sf blur422}  &0.9444(6)     &0.9549(3)     &0.9608(2)    &0.9707(1)\\
    \bottomrule[0.6pt]
    \end{tabular*}\\[2pt]
    \end{minipage}
     \centering{$\varepsilon=10^{-2}$}
      \begin{minipage}[t]{1\textwidth}
       \begin{tabular*}{\linewidth}{lp{2.0cm}p{2.0cm}p{2.0cm}p{2.0cm}p{2.0cm}p{2.0cm}}
        \toprule[0.6pt]
     	                 &$\tau=1.005$    &$\tau=1.1$    &$\tau=1.2$   &$\tau=2.0$\\ \midrule[0.6pt]
     {\sf rice}     &0.7989(10)     &0.8198(6)     &0.8267(5)    &0.8541(2)\\
     {\sf mri}      &0.8564(21)     &0.8652(15)    &0.8729(12)   &0.8985(6)    \\
     {\sf blur30}   &0.8019(56)     &0.8285(38)    &0.8564(27)   &0.9312(9)\\
     {\sf blur422}  &0.9230(21)     &0.9308(13)    &0.9367(9)    &0.9501(4)\\
          \bottomrule[0.6pt]
          \end{tabular*}\\[2pt]
          \end{minipage}
          \centering{$\varepsilon=10^{-3}$}
      \begin{minipage}[t]{1\textwidth}
       \begin{tabular*}{\linewidth}{lp{2.0cm}p{2.0cm}p{2.0cm}p{2.0cm}p{2.0cm}p{2.0cm}}
        \toprule[0.6pt]
     	                 &$\tau=1.005$    &$\tau=1.1$    &$\tau=1.2$   &$\tau=2.0$\\ \midrule[0.6pt]
           {\sf rice}   &0.7288(62)     &0.7369(46)	  &0.7430(38)	  &0.7709(18)     \\
           {\sf mri}    &0.8121(141)    &0.8179(105)  &0.8223(85)     &0.8378(41)	\\
           {\sf blur30} &0.6083(241)    &0.6183(216)  &0.6286(194)	  &0.6889(110)     \\
           {\sf blur422}&0.8861(182)    &0.8931(119)  &0.8972(92)     &0.9118(39)	\\
          \bottomrule[0.6pt]
          \end{tabular*}\\[2pt]
          \end{minipage}
\end{table}

Besides the smallest relative errors defined by \eqref{relerr},
Table \ref{tab3} also lists the relative errors
of the corresponding best regularized solutions
obtained by {\sf JBDGCV}, {\sf JBDWGCV} and {\sf JBDQR}, which are
defined by
$$
\frac{\|x_k^{reg}-x_{true}\|}{\|x_{true}\|}
$$
and marked "{\sf no $L$}" in the parentheses that follow the matrix names.
We can see that for these four problems the solution accuracy
of {\sf JBDQR} is considerably higher than
that of {\sf JBDGCV} and {\sf JBDWGCV},
no matter which relative error is used.
From the table, it is clear that the estimates for $k^*$ by
the L-curve criterion are quite rough and considerable
underestimates except for {\sf blur30} with $\varepsilon=10^{-3}$.
This indicates that the L-curve criterion does not work well
for determining $k^*$ for difficult two dimensional problems.
The fundamental cause is that $\|\bar{B}_ky_k\|$ still increases slowly
even after $k>k^*$, such that the curve of
$\left(\log(\|B_ky_k-\beta_1e_1\|),\log(\|\bar{B}_ky_k\|)\right)$
does not form a good L-shape.

Since $\|e\|$ is known for the above test problems,
we also use the discrepancy principle criterion \eqref{disc}
to estimate the optimal $k^*$.
We report the results obtained when $\tau=1.005,\
1.1, \ 1.2$  and $2.0$ in Table~\ref{addtab2}.
We can see that, for the four $\tau>1$'s, the regularization parameters
determined by the discrepancy principle have big differences for
both the solution accuracy and the estimates for $k^*$.
It is obvious that the estimates are much better for $\tau=1.005$
than those when $\tau=2$. Again, this indicates
that $\tau=2$ is definitively a very bad choice.

\begin{figure}[!htp]
\begin{minipage}{0.48\linewidth}
  \centerline{\includegraphics[width=6.0cm,height=4cm]{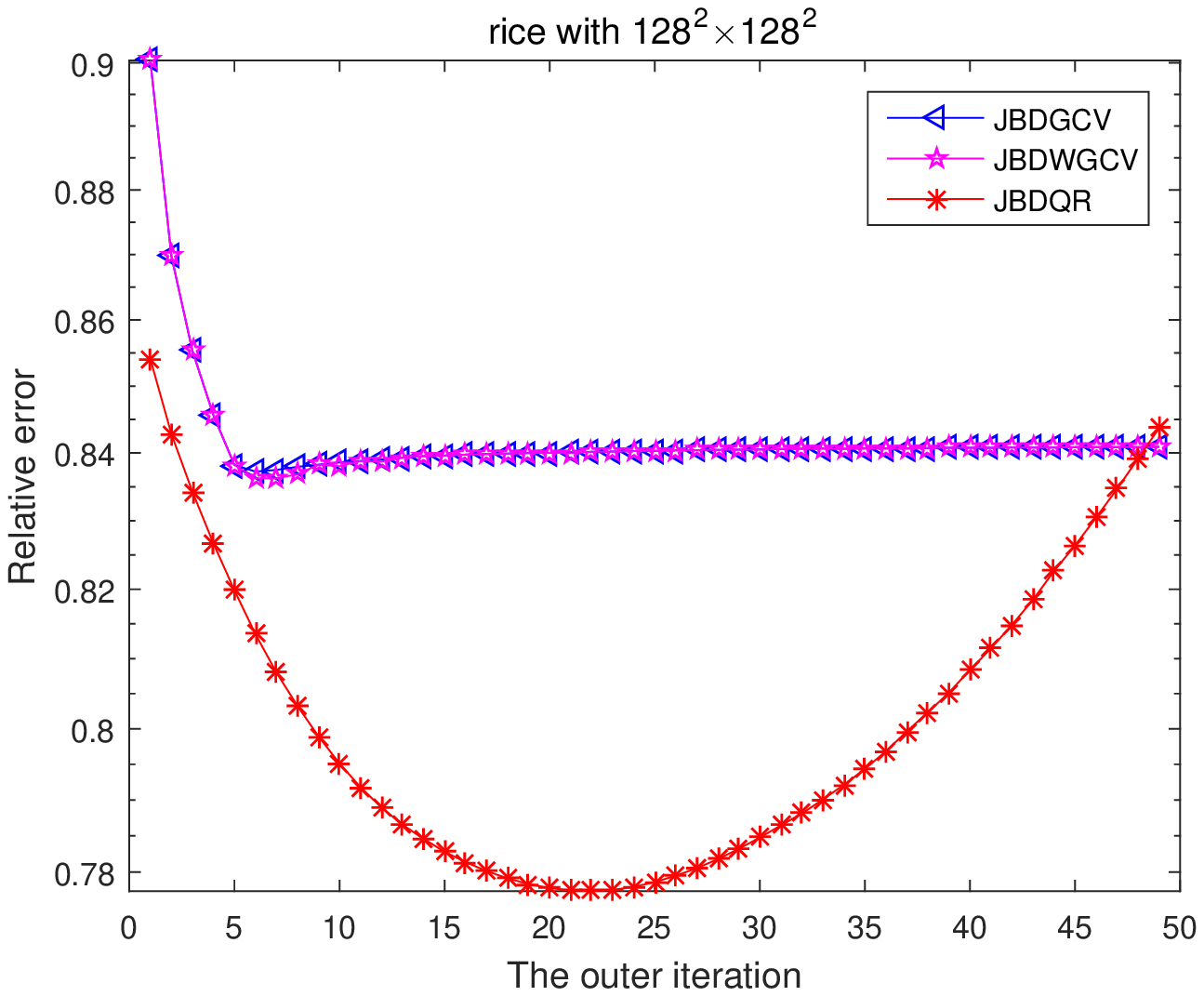}}
  \centerline{(a)}
\end{minipage}
\hfill
\begin{minipage}{0.48\linewidth}
  \centerline{\includegraphics[width=6.0cm,height=4cm]{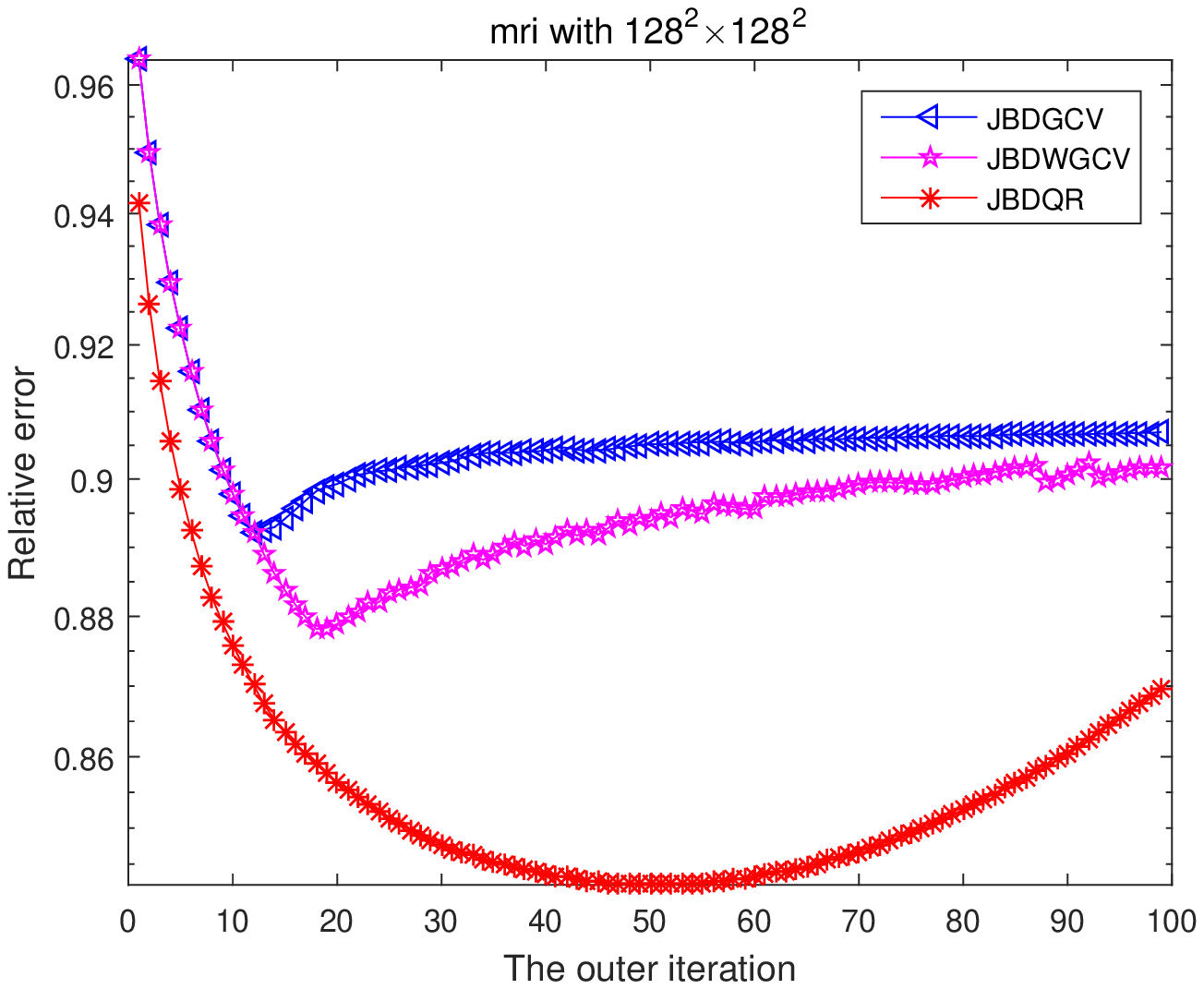}}
  \centerline{(b)}
\end{minipage}
\vfill
\begin{minipage}{0.48\linewidth}
  \centerline{\includegraphics[width=6.0cm,height=4cm]{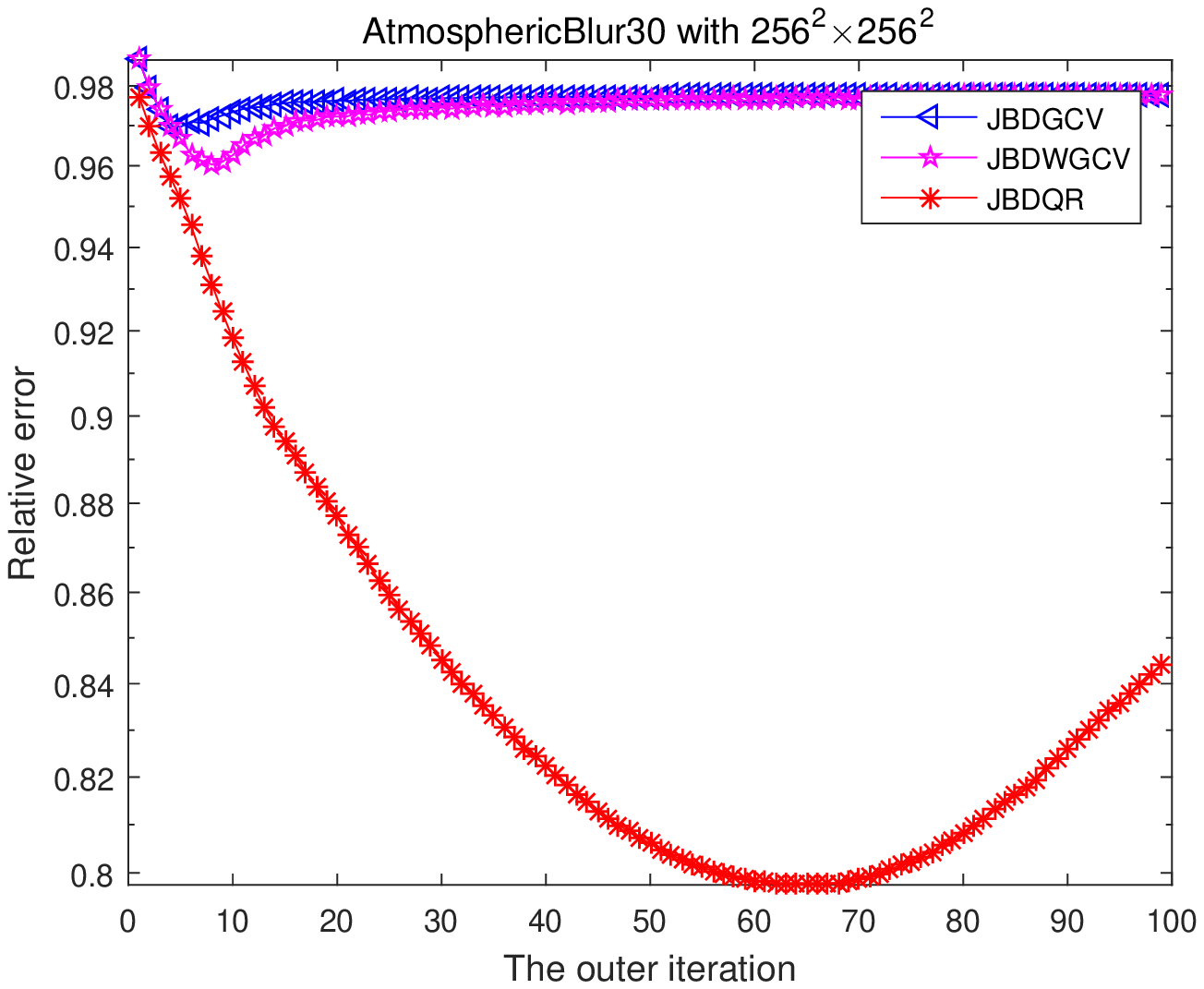}}
  \centerline{(c)}
\end{minipage}
\hfill
\begin{minipage}{0.48\linewidth}
  \centerline{\includegraphics[width=6.0cm,height=4cm]{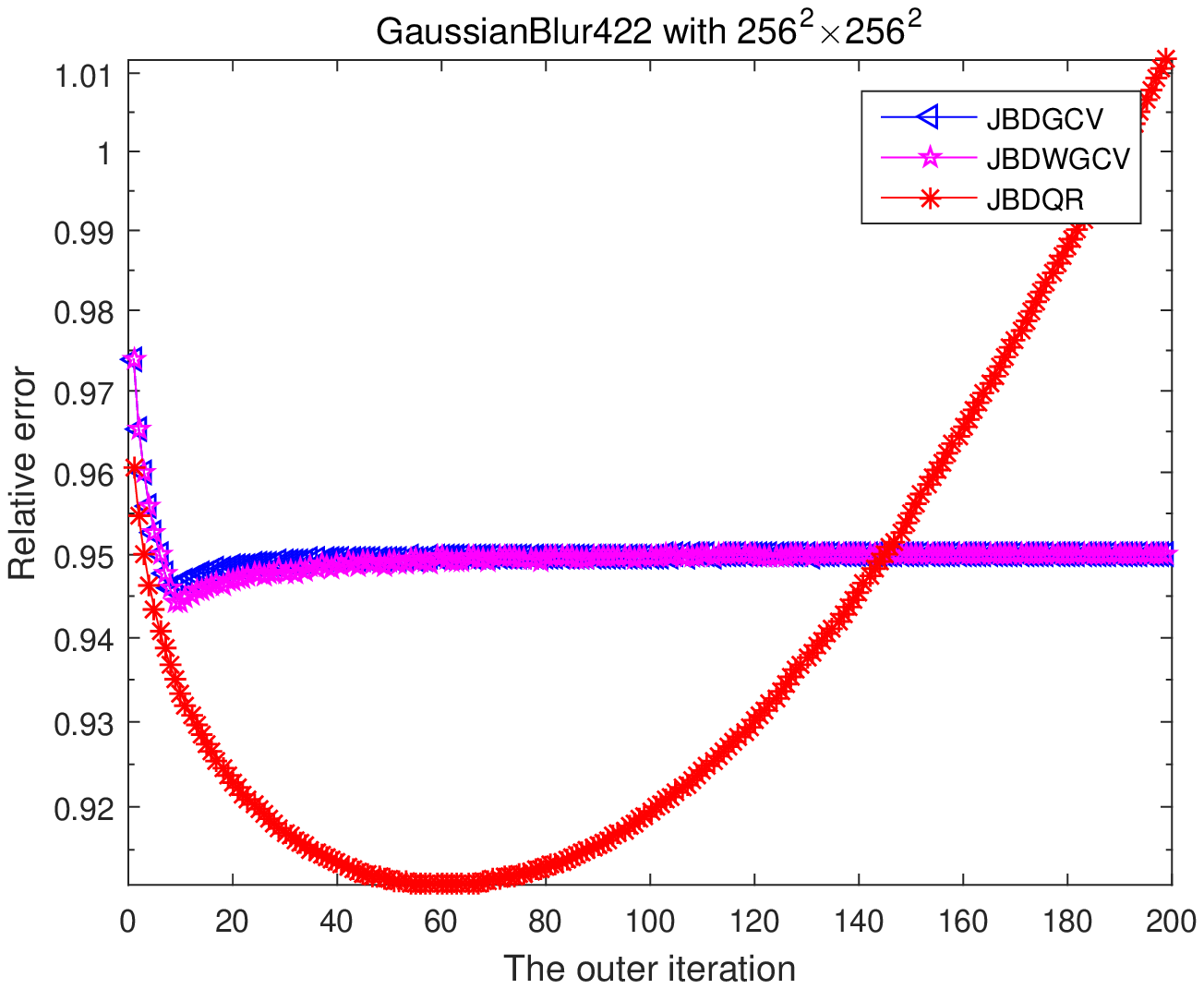}}
  \centerline{(d)}
\end{minipage}
\caption{The relative errors by {\sf JBDQR}, {\sf JBDWGCV} and {\sf JBDGCV}
with $\varepsilon=10^{-2}$:
(a) {\sf rice}; (b) {\sf mri}; (c) {\sf blur30}; (d) {\sf blur422}.}
\label{fig3}
\end{figure}

Figure \ref{fig3} draws the convergence processes of
{\sf JBDQR}, {\sf JBDGCV} and {\sf JBDWGCV}
for $\varepsilon=10^{-2}$.
We can see that the best regularized solutions by {\sf JBDQR} are more accurate
than the counterparts by {\sf JBDGCV} and {\sf JBDWGCV};
the convergence curves
of {\sf JBDGCV} and {\sf JBDWGCV} first decrease with $k$,
then increase for a while and finally stabilize, but {\sf JBDQR}
has typical semi-convergence phenomenons for all
the problems.

A final note on Table \ref{tab3} and Figure \ref{fig3} is that
the best regularized solutions by {\sf JBDWGCV} are slightly more accurate
than those by {\sf JBDGCV}, which are different from the previous results in
the one dimensional case.

\begin{figure}[!htp]
\begin{minipage}{0.48\linewidth}
  \centerline{\includegraphics[width=6.0cm,height=4cm]{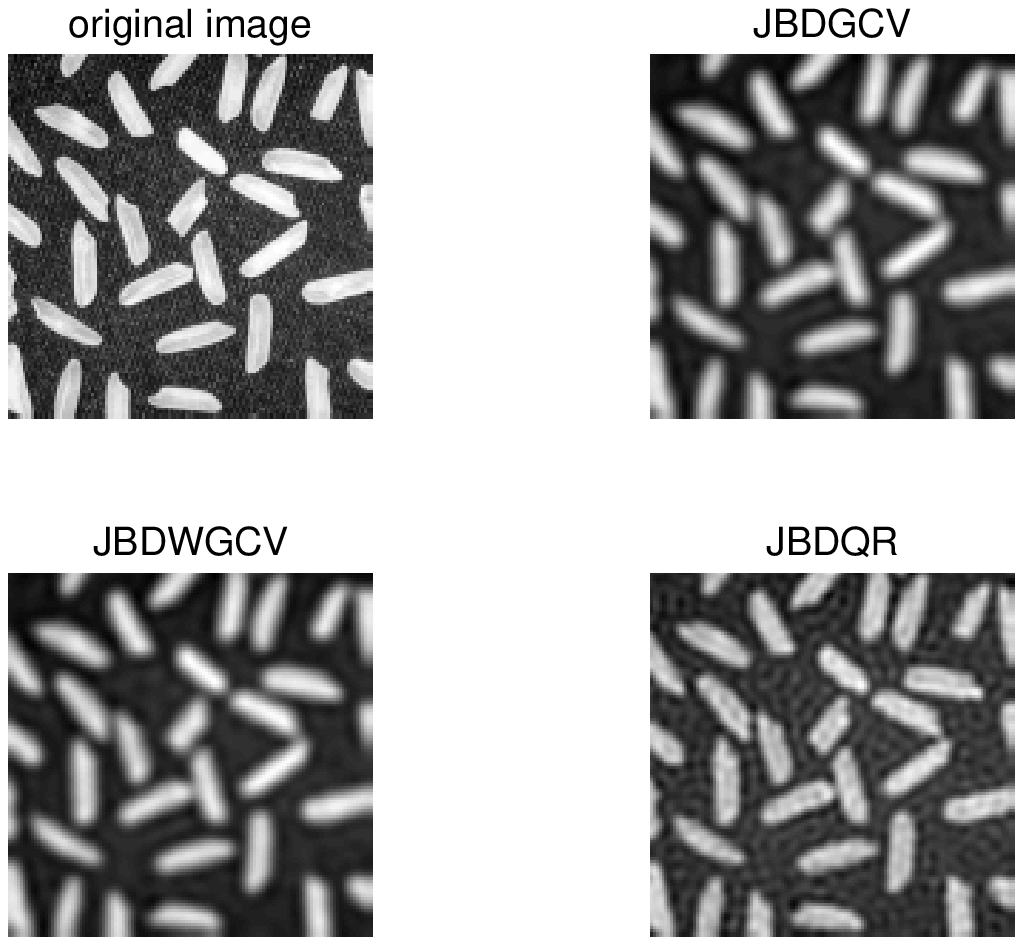}}
  \centerline{(a)}
\end{minipage}
\hfill
\begin{minipage}{0.48\linewidth}
  \centerline{\includegraphics[width=6.0cm,height=4cm]{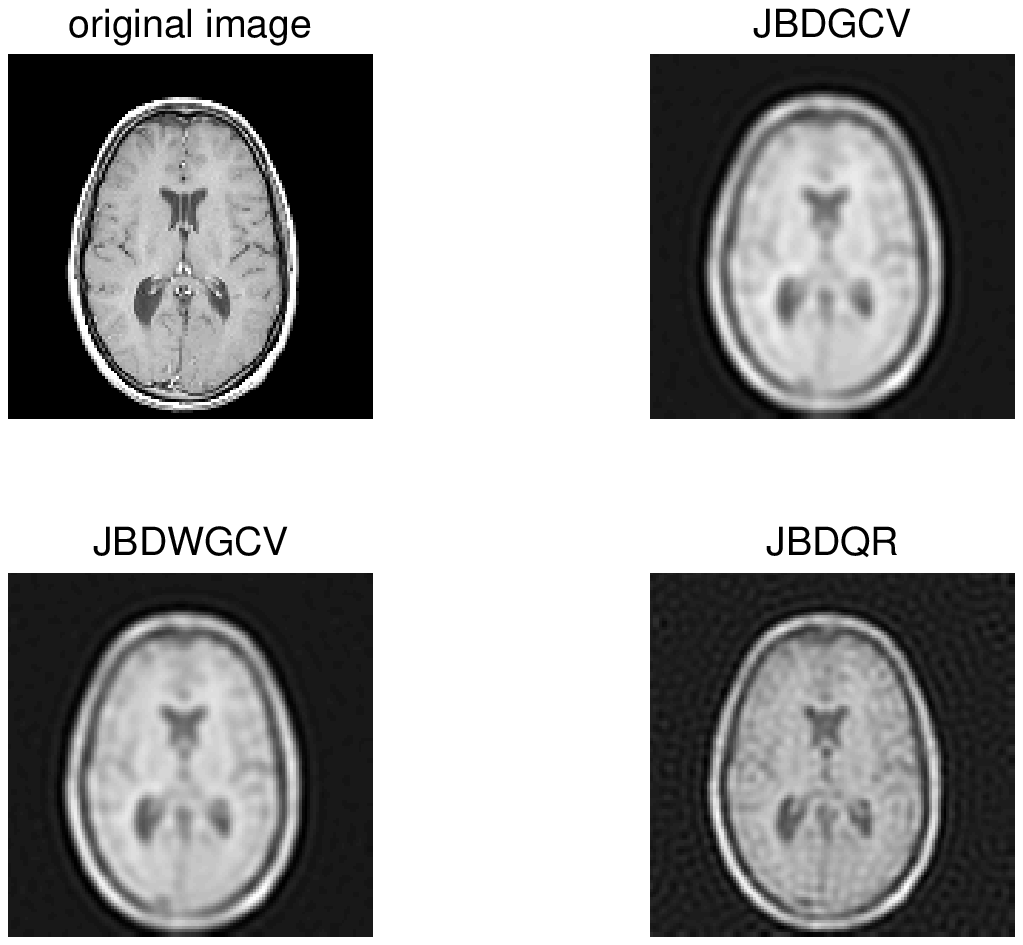}}
  \centerline{(b)}
\end{minipage}
\vfill
\begin{minipage}{0.48\linewidth}
  \centerline{\includegraphics[width=6.0cm,height=4cm]{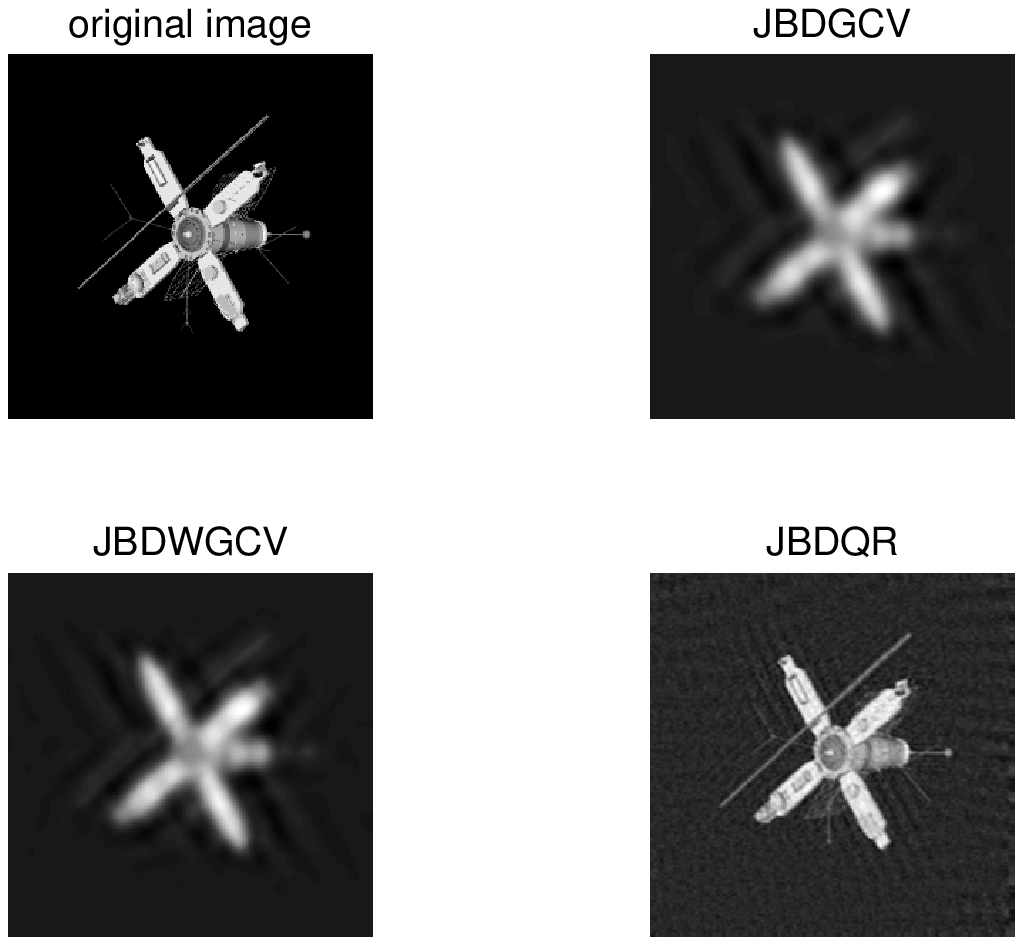}}
  \centerline{(c)}
\end{minipage}
\hfill
\begin{minipage}{0.48\linewidth}
  \centerline{\includegraphics[width=6.0cm,height=4cm]{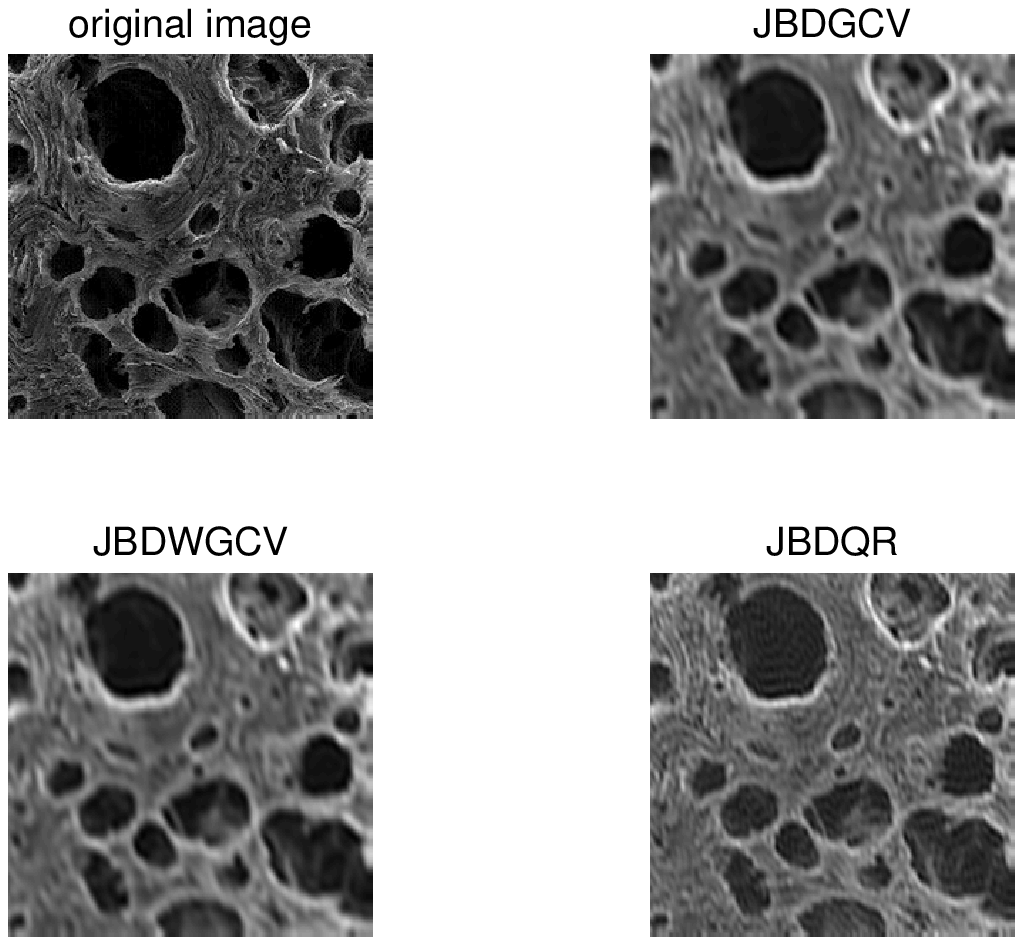}}
  \centerline{(d)}
\end{minipage}
\caption{The exact images and the reconstructed images for the four two dimensional test problems
with $\varepsilon=10^{-2}$: (a) {\sf rice}; (b) {\sf mri}; (c) {\sf blur30}; (d) {\sf blur422}.}
\label{fig7}
\end{figure}

Figure \ref{fig7} draws the exact images and
the reconstructed images for the four test problems
with $\varepsilon=10^{-2}$. Clearly, the reconstructed
images by {\sf JBDQR} are at least as sharp as those
by {\sf JBDGCV} and {\sf JBDWGCV}, and some of the former ones can be
much sharper than the latter, e.g., {\sf blur30}.

For $\varepsilon=5\cdot 10^{-2}$ and $10^{-3}$, we have similar findings
to those in Figures~\ref{fig3}--\ref{fig7}.

\section{Conclusions}

In this paper, we have proposed a joint bidiagonalization based algorithm
for solving large scale linear discrete ill-posed problems in general-form
regularization.
This algorithm is different from the hybrid projection based method proposed in
\cite{kilmer2007}, which exploits the same joint bidiagonalization process and
explicitly regularizes each projected problem
generated at every iteration.

We have analyzed the proposed algorithm and established a number of theoretical
results. Particularly, we have proved that
the iterates take the desired and attractive form of filtered GSVD expansions.
The results rigorously show that the algorithm must possesses the
semi-convergence property and get insight into the regularizing effects of
the algorithm. Our algorithm is simpler and easier to implement
than the hybrid one, and it is also more reliable and behaves regularly
than the latter.

We have made numerical experiments on a number of problems
to justify numerous aspects of the proposed algorithm,
e.g., solution accuracy and reliability. The results have illustrated
that our algorithm often computes considerably more
accurate regularized solutions than the hybrid algorithm.

There are some important unsolved problems. As we have seen, a bottleneck
of our algorithm and the hybrid one is solve a large scale least squares
problem at each outer iteration, which may be costly, especially
when the solution accuracy of these problems is high.
It is unclear if the solution accuracy can be relaxed substantially,
at least at some outer iterations, similar to
the randomized SVD algorithms proposed in
\cite{jiayang18} that solve the general-form
regularization problem \eqref{tik2}. If they could be solved with considerably
relaxed accuracy, we shall gain much, and the overall efficiency of the
algorithm can be improved substantially. The solution accuracy requirement on
the inner least squares problems will constitute our forthcoming work.

\end{document}